\documentclass[final]{article}
\usepackage{amssymb,amsmath,amsfonts,amsthm,mathrsfs,graphicx,color}
\usepackage[colorlinks=true,linkcolor=blue]{hyperref}

\usepackage{graphicx,epsf}
\newtheorem{theor}{Theorem}

\newtheorem{lem}{Lemma}[section]
\newtheorem{prop}{Proposition}
\newtheorem{remark}{Remark}

\newcommand{\R}{{\mathbb R}}

\newcommand{\calL}{\mathcal L}

\newcommand{\dist}{{\rm dist}}
\newcommand{\dive}{{\rm div}}

\newcommand{\epsl}{\varepsilon}

\newcommand{\refq}[1]{~(\ref{#1})}
\newcommand{\rmd}{\mbox{d}}

\newcommand{\dtyt}{\mbox{\bf d}(t)}

\newcommand{\sauf}{\setminus}
\newcommand{\sgn}{{\rm sign}}
\newcommand{\ssl}{\underline}
\newcommand{\supp}{\mbox{supp }}
\newcommand{\teta}{\theta}
\newcommand{\Om}{\Omega}
\newcommand{\om}{\omega}
\newcommand{\ove}{\overline}
\title
      {Recovering time-dependent inclusions in anisotropic heat conductive bodies}
\title{Uniqueness of time-dependent inclusions in anisotropic heat conductive bodies}

\author{O. Poisson\thanks{Aix Marseille Universit{\'e}, I2M, UMR CNRS 6632,
 France ({\tt  olivier.poisson@univ-amu.fr}).}}

\begin{document}
\maketitle
%

\begin{abstract}
 We consider an inverse boundary value problem for the heat equation with a nonsmooth
 coefficient of conductivity which models the displacement of a moving body inside
 a nonhomogeneous background.
 We prove the uniqueness of the moving inclusion from the knowledge
 of the Dirichlet-to-Neumann operator by using a dynamical probe method.
\end{abstract}


\medskip
{\bf AMS} : 35R30, 35K05.

\pagestyle{myheadings}
\thispagestyle{plain}

\section{Introduction}
\label{s.1}
\subsection{Inverse heat conductivity problem}
 Let $T>0$ and let $\Om$ be a bounded domain in $\R^3$, with a lipschitzian boundary
 $\Gamma = \partial\Om$.
 Let us consider the anisotropic heat equation
\begin{equation}
 \label{eq.heat1}
 \partial_t v - \dive\, ({\bf a} \nabla v) =0 \quad {\rm in} \quad \Om_{0,T}\equiv \Om\times (0,T),
\end{equation}
 where the operators $\dive$, the divergence, and $\nabla$, the gradient,
 are relative to the spatial variable $x$.
 In our model, the conductivity ${\bf a} =(a_{ij})_{1\le i,j\le 3}$ is a $3\times 3$
 real symmetric matrix with positive bounded measurable coefficients of $x$.
 It satisfies the uniform elliptic condition:\\
 there exists $\gamma_\infty>0$ such that
\begin{equation}
\label{h.gamma}
 \gamma_\infty^{-1} |\xi|^2 \le {\bf a}\xi\cdot \xi \le \gamma_\infty |\xi|^2, \quad \xi\in\R^3.
\end{equation} 
 It is well-known that, for all $f\in L^2(0,T;H^{1/2}(\Gamma))$ and $v_0\in L^2(\Om)$,
 there exists only one solution $v=v({\bf a},v_0;f)\in H^1((0,T);L^2(\Om))\cap L^2((0,T);H^1(\Om))$
 of\refq{eq.heat1} with the following initial boundary value problem:
\begin{eqnarray}
\left\{ \begin{array}{rll}
 v &=& f  \quad {\rm on} \quad \Gamma_{0,T} \equiv \Gamma\times (0,T), \\
 v\big|_{t=0} &=& v_0  \quad {\rm on} \quad \Om.
\end{array} \right.
\label{parab.v}
\end{eqnarray}
  \textcolor{black}{See for example the book of Wloka\cite{WLO.PAR}.}
 Then, we can define the Dirichlet-to-Neumann map (D-N map) as
$$
 \Lambda_{{\bf a};v_0} : L^2((0,T);H^{1/2}(\Gamma)) \ni f \mapsto {\bf a}\nabla v({\bf a},v_0;f) \cdot\,\nu
  \in  L^2((0,T);H^{-1/2}(\Gamma)),
$$
 where $\nu$ denotes the outer unit normal to $\Gamma$.
 In physical terms, $f=f(t,x)$ is the temperature distribution on the boundary and
 $\Lambda_{{\bf a},v_0}(f)$ is the resulting heat flux through the boundary.
  
\medskip
 In this article we are concerned with the Calder\'on inverse problem for\refq{eq.heat1} which is
 to determine ${\bf a}$ from the knowledge of the D-N map $\Lambda_{v_0,{\bf a}}$.
 The conductivity ${\bf a}$ consists in a non necessarily smooth background and
 an unknown inclusion $t\mapsto D_t\subset\Om$ which moves continuously inside the body $\Om$.
 Thus, in our inverse problem, the function ${\bf a}|_{\Om\sauf D_t}$ coincides with
 a measurable real matrix-function ${\bf b}\in L^\infty(\Om)$ which satisfies\refq{h.gamma}
 and represents the conductivity of a background medium, and so, is known.
 The inverse problem we address is to determine the moving inclusion
 $D=\cup_{0\le t\le T}\:  (D_t\times \{t\}) \: \subset \Om_{0,T}$ from the knowledge of
 $\Lambda_{{\bf a},v_0}$.
\begin{remark}
 In our problem the value of the conductivity inside the inclusion, ${\bf a}|_{D_t}$,
 and the initial value of $v$, $v_0$, are unknown but the article does not deal with
 their determination.
\end{remark}
\subsection{Main assumptions}
 The two following assumptions were already considered by   \textcolor{black}{several authors
 in the isotropic situation \cite{DAI.PRO},\cite{ISA.REC},\cite{POI.REC}}.
\begin{quote}
\label{as.H0}
 (H0): there exists a positive constant $\delta_1$ such that
$$
 {\rm (H0a)}: \quad  {\bf b}^{-1}-{\bf a}^{-1} \le -\delta_1<0,\:
  {\bf b}-{\bf a}\ge \delta_1>0  \quad  {\rm in} \: D, \quad
$$
 \centerline{\bf or}
$$
  {\rm (H0b)}: \quad {\bf b}^{-1}-{\bf a}^{-1}  \ge \delta_1>0,\:
  {\bf b}-{\bf a}\le -\delta_1<0 \quad  {\rm in} \: D.
$$
%

(H1): for all $t\in [0,T]$, the set $\R^3\sauf\ove{D_t}$ is connected. 
\end{quote}

 Because of technical limitations of our method when ${\bf b}$ is not sufficiently smooth,
 we need some additional geomerical assumptions on $D$.
 For a point $x\in\R^3$ and a non-empty set $E\subset\R^3$ we denote by $\rmd(x,E)$
 the quantity $\inf_{z\in E} |x-z|$ and by $|E|$ the Lebesgue-measure of $E$.
  
\begin{quote}
(H2): $t\mapsto D_t$ is lipschitzian in the following sense:\\
 there exists $K_D>0$ such that for all $x\in \ove{\Om}$ the mapping $t\mapsto \rmd(x,\Om\sauf D_t)$
 is lipschitzian in $[0,T]$ with lipschitzian constant $K_D$ and the mapping $t\mapsto \rmd(x,D_t)$
 is lipschitzian at all $s\in [0,T]$ such that $D_s\neq\emptyset$ with lipschitzian constant $K_D$.
\end{quote}
\begin{quote}
(H3):\\
(H3a): for all $t\in [0,T]$, $D_t$ satisfies the exterior cone property, i.e., \\
 there exists $\rho(t)>0$ such that for all $z\in \partial D_t$, there exists an open cylindrical cone
 $C_o(z,\rho)\subset \R^3\sauf \ove{D_t}$ with summit $z$, hightness $\rho$ and volume $\rho^3$,\\
  \centerline{\bf and}\\
(H3b):  there exists $L_D\in (0,1)$ such that\\
$|D_t\cap B(z,r)| \ge L_D \min(|D_t|,|B(z,r)|)$, $\forall r>0$, $z\in \partial D_t, \, t\in [0,T]$.
\end{quote}

\medskip

 The Runge approximation method  in the dynamical probe method is based on the uniqueness
 property (UC) which holds if the conductivity is constant.
 \textcolor{black}{However V. Isakov has shown that (UC) can fail if the conductivity is not sufficiently
 regular \cite{ISA.REC}.}
 Therefore we add the following assumption on ${\bf b}$:

\begin{quote}
 (UC) in $\Om$ - Let a sufficiently smooth domain $\om\subset\Om$,
 $a<b$ and let $u\in H^1(0,T$; $L^2(\Om))$ $\cap$ $L^2(0,T;H^1(\Om))$ such that
 $\partial_t u-\dive({\bf b}\nabla u)=0$ in $\om\times (a,b)$ and
 $u={\bf a}\nabla u \, \cdot \nu=0$ on $S\times (a,b)$, where $S$ is an non-empty open subset
 of $\partial\om$.
 Then, necessarily, $u=0$ in $\om\times (a,b)$.
\end{quote}
\begin{remark}
\label{rem.UC}
 The above definition of (UC) is independent of the choice of the time-interval $[0,T]$
 since in our work we assume that ${\bf b}$ does not depend on the variable $t$.
\end{remark}
\begin{remark}
\label{rem.UC2}
 Condition (UC) holds if ${\bf b}$ is lipschitzian or piecewise smooth:
 see \textcolor{black}{ the results of} Vessella \cite[chap 5]{VES.QUA}.
\end{remark}
%
%
\subsection{Main Result}
 Here we state our uniqueness result for the above inverse problem.
 Let $v_0$, $v'_0\in L^2(\Om)$, two conductivities ${\bf a},{\bf a}'$ satisfying
 (H0)-(H3) and (UC). Let $D'$ the inclusion related to ${\bf a}'$.

\begin{theor}
\label{t.main}
 Assume that $\Lambda_{v_0,{\bf a}}=\Lambda_{v'_0,{\bf a}'}$.
 Then, $D=D'$.
\end{theor}
\begin{remark}
 Our proof of Theorem \ref{t.main} is not completely constructive, although it is based on
\textcolor{black}{  the same dynamical method developed by the author who showed a (theoretical)
 reconstruction of $D$ from the knowledge of  $\Lambda_{v_0,{\bf a}}$ \cite{POI.REC}.}
\end{remark}

\begin{remark}
\label{rem.5}
 We shall proof Theorem \ref{t.main} with the following assumption:
$$ \ove{D(t)}\subset \Om, \quad t\in [0,T] .$$
 Therefore we replace (H1) by:

(H1'): one has $\ove{D(t)}\subset \Om$, and the set $\Om\sauf\ove{D_t}$ is connected,
 for all $t\in [0,T]$.

\medskip
 The general proof of Theorem \ref{t.main} where $D(t)$ may touch $\partial\Om$ is
 easily get from the following modification on the case (H1'):
\begin{itemize}
\item We consider a large smooth bounded domain $\Om'$ containing $\ove{\Om}$
 and we put ${\bf b}=I_3$ (the $3\times 3$ identity matrix) in $\Om'\sauf\Om$.
\item (If necessary)\footnote
{the question that (UC) in $\Om$ would imply (UC) in $\Om'$ is out of the scope of this article}
(UC) is assumed with $\Om$ replaced by $\Om'$.
\end{itemize}
\end{remark}

\begin{remark}
\label{rem.H0}
 The proof of Theorem \ref{t.main} will show that (H0) can be extended to the following situation:
\begin{quote}
 (H0') There exist positive constants $\epsl_0,\delta_1$, such that for $(x,t)\in \bar D$,
$$
   {\bf b}^{-1}(x) -{\bf a}^{-1}(x) \le -\delta_1<0,\:
 {\bf b}(x)-{\bf a}(x)\ge \delta_1>0  \quad {\rm if} \:  \rmd(x,\partial D_t)\le \epsl_0, \quad
$$
 \centerline{or}
$$
  {\bf b}^{-1}(x)-{\bf a}^{-1}(x)  \ge \delta_1>0,\:
 {\bf b}(x)-{\bf a}|_{D_t}(x)\le -\delta_1<0 \quad {\rm if} \: \rmd(x,\partial D_t)\le \epsl_0,
$$
\end{quote}
\end{remark}

\subsection{Outline}
 In Section \ref{s.2} we recall the basis of the dynamical probe method, the Runge approximation
 method and we construct indicator and pre-indicator functions from special Cauchy boundary data.
 In Section \ref{s.3} we state the lower and upper estimates on the indicator function
 from which the proof of our main Theorem \ref{t.main} can be achieved in Section \ref{s.4}.
 In Section \ref{s.5} we develop the technical results on which the proof of the estimates of
 Section \ref{s.3} is based.
 
\section{The dynamical probe method (DPM) with special solutions of the heat equation}
\label{s.2}
\subsection{Notations}
 Let us give some notations for this paper. For $E\subset\R^3$, $a<b$, and for
 $U\subset\R^3\times\R$, we put $E_{a,b}=E\times (a,b)$ and
 $U_t \equiv\{x\in \R^3 \; (x,t)\in U\}$.\\
 For non-negative integers $p,q$ or $p=1/2$, $H^p(\Om)$
 $H^p(\partial\Om)$ and $H^{p,q}(\Om_{(a,b)})$ denote the usual Sobolev spaces where
 the superscripts $p$ and $q$ indicate the regularity with respect to $x$ and $t$, respectively.
 For an open set $U \subset \R^4$ with Lipschitz boundary $\partial U$, $H^{p,q}(U)$
 is defined likewise.
 More precisely, $g\in H^{p,q}(U)$ if and only if there exists $G\in H^{p,q}(R^4)$ with
 $G = g$ in $U$. If it is the case, $\|g\|_{H^{p,q}(U)}$ is defined to be
$$ \|g\|_{H^{p,q}} := \inf \|G\|_{H^{p,q}(R^4)}, $$
 where the infimum is taken over all $G$ such that $G = g$ in $U$.
 Let $X$ be a normed space of functions. A function $f(x,t)$ is said to be in $L^2((0,T); X)$
 if $f(\cdot,t)\in X$ for almost all $t\in (0,T)$ and
$$ \|f\|^2_{L^2((0,T); X)} := \int_0^T \|f(\cdot,t)\|^2_{L^2(X)} \rmd t < \infty. $$
\textcolor{black}{ (For more details, we refer to J.L. Lions and E. Magenes \cite{LIO.NON}).}\\
 We write $\calL_{\bf a}:=\partial_t - \dive\, ({\bf a} \nabla\cdot)$, so
 $\calL_I:=\partial_t-\Delta$ for the homogeneous case.
 Similarly, we consider operator for the backward related heat equation,
 $\calL_{\bf a}^*:= -\partial_t- \dive\, ({\bf a} \nabla\cdot)$.\\
 We denote by $B(r)$ any ball of radius $r>0$ in $\R^3$.
 The open ball $\{x\in \R^3;|y-x|<r\}$, $r>0$, is denoted $B(y,r)$.\\
  We denote by $\dtyt$ the distance between $y(t)$ and $D_t$ if $D_t\neq \emptyset$,
 i.e., $\dtyt = d(y(t),D_t)$.
 If  $D_t=\emptyset$ then we put $\dtyt=+\infty$, $1/\dtyt=0$.\\
 If $\xi\in \R^3$ then $|\xi|$ denotes the euclidian norm of $\xi$ and if
 ${\bf m}$ is a $3\times 3$ real matrix then $|{\bf m}|:= \sup_{\xi\in\R^3,\; |\xi|=1}|{\bf m}\xi\cdot\xi|$.

\subsection{Brief history of the determination of an inclusion from the D-N map}
 The determination of a sufficiently smooth moving inclusion inside an homogeneous body
 was stated \textcolor{black}{by A. Elayyan and V. Isakov \cite{ELA.UNI}. Their proof is by contradiction.}
 DPM for\refq{eq.heat1} is an extension of Ikehata's probe method which was developed
 for the elliptic equation $\dive({\bf a}\nabla v)=0$ \textcolor{black}{ where ${\bf a}$
 may be tensorial \cite{IKE.SIZ}.}
 In the parabolic situation, DPM was firstly presented by Y. Daido, H. Kang and G. Nakamura
 in the case where the background is homogeneous and $D_t\in C^2$ for all $t$ \textcolor{black}{\cite{DAI.PRO}.}
 \textcolor{black}{But there, although a part of DPM works} for all spatial dimension $n$, the reconstruction of $D$
 was proved only in the case $n=1$.
 \textcolor{black}{DPM of Y. Daido, H. Kang and G. Nakamura made} the Runge approximation of
 the fundamental solution of the operator $\calL_I$. 
  \textcolor{black}{(Note that an error in this work was corrected by V. Isakov, K. Kim, G. Nakamura \cite{ISA.REC}).}
 \textcolor{black}{Unlike to the DPM of Y. Daido, H. Kang, G. Nakamura, extending the method of
 A. Elayyan and V. Isakov by a more quantitative version which requires more regularity,
 M. Di Cristo and S. Vessella proved the log-stability of $\Lambda_{{\bf a},0}\mapsto D$ in the scalar
 case (${\bf a}=aI_3$) \cite{CRI.STA}}.

 Returning to DPM, the author \textcolor{black}{used}  "special solutions"
 for the classical heat operator which are more convenient functions than the basic
 fundamental solutions $\Gamma(x-y,t-s)$, because their behaviour in time and space
 are sufficiently separated \textcolor{black}{\cite{POI.REC}}.
 Since the background  \textcolor{black}{was} homogeneous, the DPM of the author
 can reconstruct any spatially irregular inclusion as in the elliptic situation \textcolor{black}{\cite{POI.REC}}.

 However, in our situation we are limited to inclusions with some kind of lipschitzian regularitiy
 (see (H2), (H3)).
 Moreover the negative part "$-C_M\dtyt^2$" in\refq{minIinfty01} makes the reconstruction
 process unclear so the proof of Theorem \ref{t.main} is by contradiction only.
 
\subsection{Runge approximation method}
\label{ss.Runge}
 The Runge approximation method for the operator unperturbed operator $\calL_I$ with
 the homogeneous conductivity ${\bf a}=I_3$ was developed \textcolor{black}{first by Y. Daido, H. Kang, G. Nakamura,
 then by V. Isakov, K. Kim, G. Nakamura  \cite{DAI.PRO}, \cite{ISA.REC}.}
 
 Let a lipschitzian curve $\Sigma:\: [0,T] \ni t\mapsto y(t)\in \R^3\sauf \ove{D_t}$
 which does not touch $D$. 
 We extend $\Sigma$ to $t\in\R$ by putting $y(t)=y(T)$ for $t\ge T$ and $y(t)=y(0)$
 for $t\le 0$.
 Then, thanks to (H1'), there exists an open set $U \subset \Om\times\R$ containing $\ove D$
 and satisfying 
%
$$
\left\{ \begin{array}{l}
 \partial U  \mbox{ is  lipschitzian},\\
 \dist(U,\Sigma) := \inf\{|x-y|;\: x\in U, \, y\in \Sigma\}>0,\\
 \Om\sauf \ove{U_t}  \mbox{ is connected, $t\in \R$}.
\end{array} \right.
$$
%
 The Runge approximation method \textcolor{black}{works} thanks to (UC) notably, and gives the following result \textcolor{black}{\cite{DAI.PRO,ISA.REC,POI.REC}.}
 For $\tau>0$ we denote $\Sigma^\tau= \cup_{t\in\R} B(y(t),1/\tau)\times\{t\}$.
\begin{prop}
\label{p.Runge}
 Assume (H1') and (UC). Let $\Sigma$ and $U$ be as above.
 Let $u\in H^{1,0}(\Om_{(0,T)})\cap H^{0,1}(\Om_{(0,T)})$ be a solution of $\calL_{{\bf b}}u=0$
 in $\Om_{(-1,T+1)}\sauf \Sigma^\tau$.
 Then for $\tau>\inf\{r>0\, | \: \dist(U,\Sigma^r)>0\} $ there exists a sequence
 $u_j \in  H^{1,0}(\Om_{(-1,T+1)}) \cap H^{0,1}(\Om_{(-1,T+1)})$ such that
%
$$
\left\{ \begin{array}{rcl}
 \calL_{{\bf b}} u_j &=& 0 \quad {\rm in} \quad \Om_{(-1,T+1)},\\
  u_j &\to & u \quad {\rm in}\quad  H^{1,0}(U)\cap H^{0,1}(U), \\
  u_j(0) &=&  u(0) \quad {\rm in}\quad  L^2(\Om).
\end{array} \right.
$$
%
\end{prop}
\subsection{Heat Kernels}
 In many researchs devoted to inverse problems for parabolic equations, 
 the background is homogeneous, i.e, ${\bf b}=I_3$.
 In such a classical situation, the heat operator is $\partial_t-\Delta$ and its usual kernel
 $\Gamma(x,t)$ has many properties, as
\begin{enumerate}
\item It is explicit:
$$
 \Gamma(x,t) = \frac1{(4\pi t)^{3/2}} e^\frac{-x^2}{4t}, \quad t>0,\quad x\in\R^3.
$$
\item It satisfies
$$ \Gamma(x,t) \le \frac{C}{\sqrt t}|\nabla\Gamma(x,t)| ,  \quad t>0,\quad x\in\R^3,$$
 for some $C>0$.
 Hence, $\Gamma(x,t)$ is small compared  to $|\nabla\Gamma(x,t)|$ as $t\to 0$.
\item
 Thanks to the Laplace transform $\int_0^\infty \cdot \: e^{-\tau^2 t} \rmd t$ of $\partial_t-\Delta$,
 we consider similarily the elliptic operator $-\Delta +\tau^2$ with  the (large) real parameter $\tau>0$.
 Its kernel $E(x;\tau)$ is explicit too:
$$
 E(x;\tau) = \int_0^\infty  \Gamma(x,t) e^{-\tau^2 t} \rmd t = \frac{e^{-\tau|x|}}{4\pi|x|}, \quad x\in\R^3.
$$
\item
 It satisfies
$$ E(x;\tau) \le \tau |\nabla E(x;\tau)| ,\quad x\in\R^3. $$
 Hence, $E(x;\tau)$ is small compared  to $|\nabla E(x;\tau)|$ as $\tau\to \infty$, uniformly in
 all bounded set of $\R^3\sauf\{0\}$. This fact was exploited \textcolor{black}{by the author \cite{POI.REC}.}
\end{enumerate}
 Let us come back to the heat equation with a general conductivity ${\bf b}$.
 We put ${\bf b}(x)=I_3$ for $x\in\R^3\sauf \ove{\Om}$.
 
 For $y\in\R^3$, we denote by $G_y\in C(\R;L^2(\R^3))$ the fundamental solution
 of
$$  \calL_{{\bf b}}G_y = \delta_{(y,0)} , $$
 which satisfies
$$ G_y(x,t)=0 ,\quad t<0 .$$
 We have the estimate:
\begin{equation}
\label{est.Gy}
 \frac{\kappa e^{-\frac{|x-y|^2}{4\kappa^2 t}}}{t^{3/2}} \le G_y(x,t) \le 
    \frac{e^{- \frac{\kappa^2|x-y|^2}{4t}}}{\kappa t^{3/2}}, \quad x\in\R^3,\; t>0,
\end{equation}
 for some constant $\kappa=\kappa({\bf b})\in (0,1)$. 
 See \textcolor{black}{the famous results of D. G. Aronson and J. Nash \cite{ARO.BOU,NAS.CON}}.\\
 For $\tau>0$ we put the Laplace Transform of $G_y(x,t)$ as
\begin{eqnarray}
\label{def.ptau}
 p_\tau(x;y) &:=& e^{-\tau^2 t}  \int_{-\infty}^{t} e^{\tau^2 s} G_y(x,t-s)  \rmd s =
   \int_0^{\infty} e^{-\tau^2 s} G_y(x,s)  \rmd s.
\end{eqnarray}
 Let us observe that $p_\tau(\cdot;y)$ belongs to $H^1_{loc}(\R^3\sauf\{y\})$ and, thanks to\refq{est.Gy},
 satisfies
\begin{eqnarray}
\label{eq.ptau}
 (- \dive\, ( {\bf b} \nabla \cdot) + \tau^2) p_\tau(\cdot;y) = \delta_y(\cdot),\\
\label{est.ptau}
 2\sqrt\pi  \frac{\kappa^2 e^{-\frac{\tau}{\kappa} |x-y|}}{|x-y|} \le p_\tau(x;y) \le
  2\sqrt\pi \frac{e^{-\kappa\tau|x-y|}}{\kappa^2 |x-y|}, \quad x\in\R^3\sauf\{y\}.
\end{eqnarray}
%
 This is also a consequence of the works of Nash and Aronson.

\subsection{Special solutions}
\label{ss.special}
 Let us consider a lipschitzian curve $\Sigma\subset\R^3\times\R$ as in Section \ref{ss.Runge},
 and fix $\teta\in (0,T)$.
 Let another positive parameter $\mu\ge 1$ that we shall precise later.

 \textcolor{black}{The author considered special solutions related to the following functions (with other notations
 and with ${\bf b}\equiv I_3$)}:
\begin{eqnarray*}
 U_{OP}(x,t) := e^{\tau^2 (T+t)} \int_0^\infty e^{\tau\mu(|t-\teta-s|-|t-\teta|)} 
  \Gamma(x-y(t-s),s) e^{-\tau^2 s} \rmd s,\\
 U_{OP}^*(x,t) := e^{-\tau^2 (T+t)} \int_0^\infty e^{\tau\mu(|t-\teta+s|-|t-\teta|)} 
  \Gamma(x-y(t+s),s) e^{-\tau^2 s} \rmd s,
\end{eqnarray*}
 \textcolor{black}{\cite{POI.REC}.}
 In fact, $U_{OP}$ and $U_{OP}^*$ are respectively solutions of the following forward
 and backward heat equations:
\begin{eqnarray*}
  \calL_I U_{OP} (x,t) &=& e^{\tau^2 (t+T)} e^{-\tau\mu |t-\teta|} p_\tau(x;y(t)) \:
   \quad {\rm in} \: \R^3\times\R,\\
  \calL_I^* U_{OP}^* (x,t) &=& e^{-\tau^2 (t+T)} e^{-\tau\mu |t-\teta|} p_\tau(x;y(t)) \: 
   \quad {\rm in} \: \R^3\times\R.
\end{eqnarray*}
 Moreover they satisfies
\begin{eqnarray*}
 U_{OP} (x,t) = \varphi(x,t) e^{\tau^2 (t+T)} e^{-\tau\mu |t-\teta|} p_\tau(x;y(t)),\\
 U_{OP}^* (x,t) = \varphi^*(x,t) e^{-\tau^2 (t+T)} e^{-\tau\mu |t-\teta|} p_\tau(x;y(t)),
\end{eqnarray*}
 such that, for some $C=C(R,\mu)>0$ and all $\tau\ge C$,
\begin{eqnarray}
\label{est1.varphi}
 \frac1C \le |\varphi(x,t)| +|\varphi^*(x,t)| \le C \quad {\rm in} \: B(0,R)\times\R ,\\
\label{est2.varphi}
 |\nabla\varphi(x,t)| + |\nabla\varphi^*(x,t)| \le C \quad {\rm in} \: B(0,R)\times\R,
\end{eqnarray}
 \textcolor{black}{\cite[Lemma 1]{POI.REC}}.
 With the general conductivity ${\bf b}$, we construct here special solutions
 $u_\tau$ and $u_\tau^*$ as follows.
 Let us put
\begin{equation}
\label{def.mtau}
 m_\tau(x,t) = M_0(\tau|x-y(t)|), \quad t\in\R,
\end{equation}
 where $M_0$ is defined by $M_0(r)=|1-r| \; 1_{|r|\le 1}$.
 Hence $m_\tau$ is a lipschitzian function with support closed to $\Sigma$ as $\tau>>1$.
 We then put, for $(x,t)\in \R^3\times\R$,
\begin{eqnarray}
\label{def.utau}
 u_\tau(x,t) &=&  \int_{s\in\R} \int_{y\in\R^3} e^{\tau^2 (s+T)} e^{-\tau\mu |s-\teta|} m(y,s)G_{y}(x,t-s) \rmd y \rmd s \\
\nonumber
  &=&  \int_{s=0}^\infty \int_{y\in\R^3} e^{\tau^2 (T+t-s)} e^{-\tau\mu |t-\teta-s|} m(y,t-s)G_{y}(x,s) \rmd y  \rmd s,\\
\label{def.utau*}
 u_\tau^*(x,t) &=&  \int_{s\in\R} \int_{y\in\R^3} e^{-\tau^2 (T+s)} e^{-\tau\mu |s-\teta|} m_\tau(y,s)G_{y}(x,s-t)\rmd y  \rmd s \\
\nonumber
  &=&  \int_{s=0}^\infty \int_{y\in\R^3} e^{-\tau^2 (T+t+s)} e^{-\tau\mu |t-\teta+s|}
   m_\tau(y,t+s)G_{y}(x,s) \rmd y  \rmd s.
\end{eqnarray}
 The functions $u_\tau$ and $u_\tau^*(x,t)$ are positive and satisfy
\begin{eqnarray}
\label{eq.utau}
 \calL_{{\bf b}} u_\tau (x,t)= e^{\tau^2 (t+T)} e^{-\tau\mu |t-\teta|} m(x,t) \: 
   \quad {\rm in} \: \R^3\times\R,\\
\nonumber
 \calL_{{\bf b}}^* u_\tau^ * (x,t)=  \: e^{-\tau^2 (T+t)} e^{-\tau\mu |t-\teta|}  m_\tau(x,t) 
  \quad {\rm in} \: \R^3\times\R.
\end{eqnarray}
\begin{remark}
\label{rem.mtau}
 If $m_\tau(x,t)$ was replaced by $\delta(x-y(t))$ then it would be difficult to make the estimation
 of $\dot y(s)\nabla_y G_{y(t-s)}(x,t)$ that would appear in the expression of $\partial_t u_\tau$.
\end{remark}
 We then expect that
%
\begin{eqnarray}
\label{asym2.utau}
  u_\tau(x,t) \stackrel{\tau\to \infty}\simeq  e^{\tau^2 (T+t)} e^{-\tau\mu |t-\teta|}  \tau^{-3} p_\tau(x,y(t)),\\
\label{asym2.utau*}
  u_\tau^*(x,t) \stackrel{\tau\to \infty}\simeq  e^{-\tau^2 (T+t)} e^{-\tau\mu |t-\teta|}   \tau^{-3} p_\tau(x,y(t)),
\end{eqnarray}
 where the meaning of "$\simeq$" will be clarified shortly.
 Since the comparison requires the time-derivatives of $u_\tau(x,t)$ or $u_\tau^*(x,t)$
 and remembering Remark \ref{rem.mtau}, we introduce the following smooth approximation
 of $p_{\tau}(x;y(t))$:
\begin{equation}
\label{def.PM}
 P_\tau(x,t) 
  :=  \int_0^{\infty} \int_{\R^3} e^{-\tau^2 s} m_\tau(y,t)  G_y(x,s)  \rmd y \rmd s.
\end{equation}
 We then put
\begin{eqnarray}
\label{def.qtau}
 q_\tau(x,t) &:=&  e^{-\tau^2 (T+t)} u_\tau(x,t) -  e^{-\tau\mu |t-\teta|} P_\tau(x,t),  \\
\label{def.qtau*}
 q_\tau^*(x,t) &:=&  e^{\tau^2 (T+t)} u_\tau^*(x,t) -  e^{-\tau\mu |t-\teta|} P_\tau(x,t).
\end{eqnarray}
 The main difficulty in the proof of Theorem \ref{t.main} is to prove that the quantity
\begin{equation}
\label{def.R0}
  R_0 :=  \int_{D} ( |\nabla q_\tau(x,t) |^2 + |\nabla q_\tau^*(x,t) |^2) \,  \rmd x \rmd t ,
\end{equation}
 is negligible compared to  $\int_{D} \tau^{-6} e^{-2\tau\mu |t-\teta|} |\nabla p_\tau(x,t)|^2 \, \rmd x \rmd t$
 or, in an equivalent way (see Lemmas \ref{lem.Caccioppoli-PM} and \ref{lem.CompPMptau}),
 to $\int_{D} \tau^{-4} e^{-2\tau\mu |t-\teta|} |p_\tau(x,t)|^2 \, \rmd x \rmd t$.
 We shall prove in Appendix the following Lemma.
\begin{lem} (Estimate of $\nabla q_\tau$ in $D_t$).
\label{lem2.maj-dqtau}
 Let $M>0$ and assume that $|\dot y|_\infty\le M$.
 There exist positive constants $C_M$, $\tau_0(\Sigma)$ such that
 if $t\in[0,T]$ and $\tau> \tau_0$ then
\begin{equation}
\label{est2.dqtau-Dt}
 \int_{D_t}  (|\nabla q_\tau(x,t)|^2+|\nabla q_\tau^*(x,t)|^2) \rmd x \le
  C_M \dtyt^2 \tau^{-4} e^{-2\tau\mu |t-\teta|}   \int_{D_t}  |p_\tau(x,y(t))|^2 \rmd x.
\end{equation}
\end{lem}
 (Remember that $\dtyt = d(y(t),D_t)$.)
 So $R_0$ is effectively "negligible" when the curve $\Sigma$ is sufficiently close to $D$
 at least at time $\teta$.
 This constraint is new compared to consequences of \refq{est1.varphi} and \refq{est2.varphi}
 (for which Assumption (H3) is in addition superfluous) and makes a theoritical reconstruction
 of $D$ problematic, as opposite to the possible reconstruction proposed \textcolor{black}{
 by the author \cite{POI.REC}}.

\subsection{Pre-indicator sequence and indicator function}
\label{ss.indicator}
 As in section \ref{ss.Runge}, we can consider sequences $(u_j)_j$ and $(u_j^*)_j$
 such that $u_j\to u_\tau$ and $u_j^*\to u_\tau^*$ in the sense of Proposition \ref{p.Runge}.
 Considering $v_j=v({\bf a},v_0;u_j|_{\Gamma_{(0,T)}})$ and the solution
 $v_\tau \in H^1((0,T);L^2(\Om))\cap L^2((0,T);H^1(\Om))$ of
\begin{eqnarray}
\left\{ \begin{array}{rll}
 \calL_{\bf a} v_\tau &=&  \calL_{{\bf b}} u_\tau,\\
 v_\tau &=& u_\tau \quad {\rm on} \quad \Gamma_{0,T}, \\
 v_\tau\big|_{t=0} &=& v_0  \quad {\rm on} \quad \Om,
\end{array} \right.
\label{parab.vtau}
\end{eqnarray}
 we put
\begin{equation}
\label{def.wtau} 
 w_\tau= v_\tau- u_\tau
\end{equation}
 and
\begin{eqnarray}
\nonumber
 I_j(\tau) &:=& \int_{\Gamma\times [0,T]} (\Lambda_{{\bf a};v_0} (u_j|_{\Gamma_{[0,T]}})
 - {\bf b}\nabla u_j \, \cdot \nu)\; u_j^*|_{\Gamma\times [0,T]} \, \rmd\sigma(x) \rmd t,\\
\label{d.Iinfty}
 I_\infty(\tau) &:=&  \int_{\Om\times [0,T]}({\bf a}-{\bf b}) \nabla v_\tau
  \nabla u_\tau^* \,  \rmd x \rmd t + \int_\Om \left[w_\tau  u_\tau^* \right]^T_0  \rmd x,
\end{eqnarray}
 where $\rmd\sigma(x)$ is the usual measure on the boundary $\Gamma$.
 The knowledge of $\Lambda_{{\bf a};v_0}$ involves that of $I_j(\tau)$'s.
 Furthermore, as for the proofs in similar situations\textcolor{black}{,} Proposition \ref{p.Runge}
 implies that
\begin{equation}
\label{lim1.Ij}
 I_j(\tau) \to I_\infty(\tau) \in\R \quad {\rm as} \quad j\to\infty,
\end{equation}
  \textcolor{black}{For a proof, see the works based on DPM \cite{DAI.PRO,POI.REC}.}
 Hence, if (UC) holds, then the knowledge of $\Lambda_{{\bf a};v_0}$ involves
 that of $I_\infty(\tau)$'s.
 
\section{Estimates on the indicator function}
\label{s.3}
 In the following results the positive constants $c,C,C_1$ may depend on $T$, $\Om$, $\mu$,
 on the the conductivity, but not on $\tau$. We indicate when they depend on an upper bound
 $M$ of $|\dot y|_\infty$ or on the initial data $v_0$. 
\begin{lem}
\label{est-lem.Iinfty}
 Under assumption (H0b) we have
\begin{eqnarray}
\label{majIinfty0}
 I_\infty(\tau)   &\le&  C \int_{D}  e^{-2\tau\mu |t-\teta|} |\nabla P_\tau(x,t)|^2 \, \rmd x \rmd t 
  \\
 \nonumber &&  + C \int_{D}   (|\nabla q_\tau |^2+ |\nabla q_\tau^* |^2) \, \rmd x \rmd t +
  10(\|v_0\|^2_{L^2(\Om)}+d_\Om)  e^{-\tau\mu \min(T-\teta,\teta)},
\end{eqnarray}
 and
\begin{eqnarray}
\label{minIinfty0}
    I_\infty(\tau)   &\ge& \frac1C \int_{D}  e^{-2\tau\mu |t-\teta|} |\nabla P_\tau(x,t)|^2 \, \rmd x \rmd t  \\
 \nonumber &&   -C \int_{D} (|\nabla q_\tau |^2+ |\nabla q_\tau^* |^2) \, \rmd x \rmd t
  - 10(\|v_0\|^2_{L^2(\Om)}+d_\Om) e^{-\tau\mu \min(T-\teta,\teta)},
\end{eqnarray}
 for some $C\ge 1$, for all $\tau>\mu+1$, $\mu>0$.
 \end{lem}
%
%
 Proof in Appendix.

\medskip
 We put also 
\begin{eqnarray}
\label{def.diamSigma}
 d_\Om &:=& \sup\{|x-y|;\: x,y\in \Om\},\\
\label{def.epslSigma}
\epsl_\Sigma &:=& \inf_{t\in [0,T]} \dtyt >0.
\end{eqnarray}

\begin{lem}
\label{lem.R9}
 Let $M>0$ and assume that $|\dot y|_\infty\le M$.
 Then, under assumption (H0b), there exist positive constants $c=c(M)$, $C_1=C_1(v_0)$,
 $C_M$, $\tau_0=\tau_0(\Sigma)$, such that
 if $\tau>\tau_0$ we then have
\begin{eqnarray}
\label{majIinfty01}
 I_\infty(\tau)   &\le& c \tau^{-4} \int_{D} e^{-2\tau\mu |t-\teta|} |p_\tau(x,y(t))|^2 \, \rmd x \rmd t \\
&& \nonumber
  + C_1 e^{-\tau\mu \min(T-\teta,\teta)},
\end{eqnarray}
 and 
\begin{eqnarray}
\label{minIinfty01}
  I_\infty(\tau)   &\ge& \frac1c \tau^{-4}  \int_{D} \big(1- C_M\dtyt^2\big) e^{-2\tau\mu |t-\teta|}
   |p_\tau(x,y(t))|^2 \, \rmd x \rmd t \\
&& \nonumber
  - C_1e^{-\tau\mu \min(T-\teta,\teta)}.
\end{eqnarray}
\end{lem}
%
%
 The proof of Lemma \ref{lem.R9} requires the developments of Section \ref{s.5}.
 Let us extend the function $\ln$ to $(-\infty,0]$ by putting $\ln(I)=-\infty$ if $I\le 0$.
\begin{lem}
\label{lem.est1-lnIinfty}
 We assume that (H0b) is true.
 Let $\teta\in (0,T)$.
 Let us fix $\mu\ge 4\kappa^{-1} d_\Om\max((T-\teta)^{-1},\teta^{-1})$.
 
 1)  Let a lipschitzian curve $\Sigma$ such that
   $\epsl_\Sigma>0$.
 We then have
\begin{eqnarray}
\label{est1.lnIinfty}
  \limsup_{\tau\to\infty} \tau^{-1} \ln(I_\infty(\tau)) \le -2\kappa \epsl_\Sigma.
\end{eqnarray}

2) Assume that $D_\teta\neq\emptyset$. 
 Let $M>0,\alpha>0$.
 Let a familly of curves as in Section \ref{ss.Runge}, $(\Sigma=\Sigma(\epsl))_{0<\epsl\le \alpha^2}$, 
 such that we have
\begin{equation}
\label{cond1.Sigma}
\left\{ \begin{array}{rcl}
 |\dot y(\cdot)|_\infty &\le& M,\\
 \dtyt &\le& \frac2{\alpha^3}\epsl \quad {\rm for}\:  |t-\teta|\le  \epsl,\\
 \dtyt &\in& [\frac1{2\alpha}|t-\teta|,\frac2{\alpha^3}|t-\teta|] 
  \quad {\rm for}\:     \epsl\le   |t-\teta|\le\alpha^2,\\
 \dtyt &\ge& \alpha/2 \quad {\rm for}\:   |t-\teta|\ge\alpha^2,
\end{array}\right.
\end{equation}
 where $\dtyt:=\rmd (y(t),D_t)$.
 Then there exists $\epsl_1 \in (0,\alpha^2]$ such that for $0<\epsl\le \epsl_1$
 we have
\begin{eqnarray}
\label{est2.lnIinfty}
  \liminf_{\tau\to\infty} \tau^{-1} \ln(I_\infty(\tau)) \ge -(8\kappa^{-1}\alpha^{-3}+4\mu)\epsl.
\end{eqnarray}
\end{lem}
 Proof in Appendix.

\section{Proof of Theorem \ref{t.main}}
\label{s.4}
 We may assume that (H0b) holds, since the case where (H0a) holds is similar.
 Thanks to Remark \ref{rem.5} we have $\ove{D_t\cup D'_t}\subset \Om$, $t\in [0,T]$.
 Let us assume that $D\neq D'$.
 Then there exists $(z,\teta)\in \Om\times [0,T]$ with $D_\teta\neq \emptyset$,
 $z\in \partial D_\teta$ and $z\not\in \ove{D'_\teta}$ or with $D'_\teta\neq \emptyset$,
 $z\in \partial D'_\teta$ and $z\not\in \ove{D_\teta}$.
 Thus, we consider for simplicity that $z\in \partial D_\teta$ and $z\not\in \ove{D'_\teta}$.
 Thanks to (H2), $t\mapsto (D_t,D'_t)$ is continuous so we consider also that $0<\teta<T$.
 In fact let us explain why can consider also that $z\in \partial D_\teta\sauf \ove{D'_t}$
 if $|t-\teta|<\beta$ for some $\beta>0$.
 If $D'_t$ is void for $|t-\teta|$ sufficiently small then it is immediate, but
 if $D'_t$ is not void for $|t-\teta|$ sufficiently small then we can't be sure that
 $\rmd(z,D'_t)>0$ when $t\simeq \teta$. However, in such a case, thanks to (H2),
 there exists a sequence $\teta_n\to \teta$ satisfying $D'_{\teta_n}\neq \emptyset$
 and $D_{\teta_n}\sauf \ove{D'_{\teta_n}}\neq \emptyset$.
 We then replace $(z,\teta)$ by another couple $(z_n,\teta_n)$ with
 $z_n\in \partial D_{\teta_n}\sauf \ove{D'_{\teta_n}}$. 
 Then, since $D'_{\teta_n}\neq \emptyset$ and thanks to (H2), we have
 $z_n \not\in  \partial D'_t$ if $t\simeq \teta_n$.\\
 So we can consider that
\begin{equation}
\label{cond0.z-teta}
 z\in \partial D_\teta\sauf \ove{D'_t} \:  \mbox{ if $|t-\teta|\le\beta$ for some $\beta>0$}.
\end{equation}
 Let us construct a familly of curves $\Sigma=\Sigma(\epsl)$ for $0<\epsl\le \alpha^2$,
 for some positive $\alpha$ such that\refq{cond1.Sigma} and
\begin{equation}
\label{cond2.Sigma}
 \epsl'_\Sigma:=\inf_{0\le t\le T} \rmd (y(t),D'_t) \ge \alpha/2
\end{equation}
 hold.
 In fact, since $D_\teta$ and $D'_\teta$ satisfy Assumptions (H1') and (H3b) then there exists
 a lipschitzian curve $\tilde y:\: [0,1]\ni s \mapsto \tilde y(s)\in \R^3$ with Lipschitz constant
 $\tilde M$ such that $\tilde y(0)=z$, $\tilde y(1)\not\in \Om$,
 $\tilde y(s)\not\in \ove{D_\teta}$ for $s\neq 0$  and $\tilde y(s)\not\in \ove{D'_t}$
 for $s\in [0,1]$ and $|t-\teta|\le \beta$.
 Thanks to (H2), (H3a) and to\refq{cond0.z-teta} we have for all $s\in [0,1]$
\begin{eqnarray}
\nonumber
 \rmd (\tilde y(s),D'_t) &\ge& \alpha- K_{D'}|t-\teta|,\\
\label{ine2.tildey}
 \rmd (\tilde y(s),D_t) &\in& [\alpha s - K_D|t-\teta|,\frac1{\alpha}s +  K_D|t-\teta|],
\end{eqnarray}
 where $\alpha>0$ is sufficiently small. We may consider that
\begin{equation}
\label{cond1.alpha}
\alpha\le\min(1,(2K_D)^{-1},(2K_{D'})^{-1},\rmd(\partial\Om,D_t),\rmd(\partial\Om,D'_t)),
 \quad t\in [0,T].
\end{equation}
 Then we have
\begin{eqnarray}
\label{ine1.tildey}
 \rmd (\tilde y(s),D'_t) &\ge& \alpha/2 \quad {\rm for} \quad  |t-\teta|\le \frac{\alpha}{2K_{D'}}, \: s\in [0,1].
\end{eqnarray}
%
 We put $y_0(t)=\tilde y(|t-\teta|/\alpha^2)$ for $|t-\teta|\le \alpha^2$ and
 $y_0(t)=\tilde y(1)$ for $|t-\teta|\ge \alpha^2$.
 From\refq{cond1.alpha},\refq{ine1.tildey} or \refq{ine2.tildey}, and since $y_0(r)\not\in \Om$
 for $|r-\teta|\ge \alpha^2$ we obtain
\begin{eqnarray}
\label{ine1.y0}
 \rmd (y_0(r),D'_t) &\ge& \alpha/2 , \quad {\rm for} \quad t,r\in [0,T],\\
 \label{ine2.y0}
\rmd (y_0(r),D_t) &\ge& \alpha/2 , \quad {\rm for} \quad t,r\in [0,T],\: |r-\teta|\ge \alpha^2.
\end{eqnarray}
 Then for all $\epsl\in (0,\alpha^2]$ we put
$$
 y(t) = \left\{ \begin{array}{rcl} y_0(\teta+\epsl) & {\rm for} &  |t-\teta|\le \epsl \\
  y_0(t)  & {\rm for} & |t-\teta|\ge \epsl.
  \end{array}
 \right.
$$
 Thanks to\refq{cond1.alpha},\refq{ine2.tildey}\refq{ine1.y0},\refq{ine2.y0} we then obtain
 all the conditions of\refq{cond1.Sigma} with $M=\tilde M/\alpha^2$, and\refq{cond2.Sigma}.
 
\medskip
 Let us denote by $I'_\infty(\tau)$ the indicator function for the conductivity ${\bf a}'$.
 Thanks to\refq{est1.lnIinfty} of Lemma \ref{lem.est1-lnIinfty} we have
$$ \limsup_{\tau\to\infty} \tau^{-1} \ln(I'_\infty(\tau)) \le -\kappa \alpha,$$
 and there exists $\epsl_1\in (0, \alpha^2)$ such that for $\epsl\in (0,\epsl_1]$ we have,
 from\refq{est2.lnIinfty},

$$   \liminf_{\tau\to\infty} \tau^{-1} \ln(I_\infty(\tau))\ge -(8\kappa^{-1}\alpha^{-3}+4\mu)\epsl. $$
 Then, $I'_\infty(\tau)\neq I_\infty(\tau)$ for all $\epsl < \min(\epsl_1,\frac{\kappa\alpha}{8\kappa^{-1}\alpha^{-3}+4\mu})$
 and $\tau$ sufficiently large.
 The result at \S \ref{ss.indicator} implies that $\Lambda_{v_0,{\bf a}}\neq\Lambda_{v'_0,{\bf a}'}$.
 \qed

\section{Technical Results}
\label{s.5}
\subsection{Basic estimates}
\begin{lem}
\label{lem.geometrie}
 Let $t\in [0,T]$ such that $D_t\neq\emptyset$.
 Then there exists a non empty finite familly $I$, and points $x_i\in D_t$,
 $i\in I$, such that
$$
 \cup_{i\in I} B_i(1/\tau)  \subset  \ove{D_t} \subset \cup_{i\in I} B_i(3/\tau),
$$
 and $B_i(1/\tau) \cap B_j(1/\tau) =\emptyset$ if $i,j\in I$, $i\neq j$, where $B_i(R)$
 denotes the open euclidian ball of radius $R>0$ and centered at $x_i$.
\end{lem}
 Proof. The lemma is a straightforwardly consequence of the compactness of $\ove{D_t}$
 and  Vitali's lemma.
 
\medskip
 We have the following proposition:
\begin{prop} 
\label{prop.Harnack}
 (Parabolic Harnack's inequality).
 There exists $c>0$ such that if $r>0$, $t\in\R$, and if $y\in\R^3\sauf B(2r)$ or
 $0\not\in (t-r^2,t+r^2)$ then we have
\begin{equation}
\label{est.Harnack-parab}
 \max_{x\in \ove{B(r)},s\in [t-\frac34r^2,t-\frac14r^2]} G_y(x,s)  \le
  c \min_{x\in \ove{B(r)},s\in [t+\frac14r^2,t+r^2]} G_y(x,s).
\end{equation}
\end{prop}
 \textcolor{black}{ For a proof, see for example the work of E.B. Fabes and D.W. Stroock \cite{FAB.NEW}.}

 Let us remember that $p_\tau$ is defined  by\refq{def.ptau}.
 From Proposition \ref{prop.Harnack}, we prove the following Lemma.
\begin{lem}
 (Elliptic Harnack's inequality).
 Let $\beta>0$. There exists $c>0$ such that for all $\tau>0$, for all ball $B(\beta/\tau)\subset\R^N$,
 if $y\not\in B(2\beta/\tau)$ we then have
\begin{equation}
\label{est.Harnack-ellip}
 \max_{x\in \ove{B(\beta/\tau)}} p_\tau(x;y)  \le c\min_{x\in \ove{B(\beta/\tau)}} p_\tau(x;y).
\end{equation}
\end{lem}
 Proof.
 Applying\refq{est.Harnack-parab} with $s=t$, $r=\beta/\tau$, we have,
 for all $x,z\in \ove{B(\beta/\tau)}$,
\begin{eqnarray*}
 p_\tau(z;y) &=& \int_0^\infty e^{-\tau^2 s} G_y(z,s) \rmd s\\
 &=& \int_{\frac12 \beta^2/\tau^2}^\infty e^{-\tau^2 (s-\frac12 \beta^2/\tau^2)} G_y(z,s-\frac12 \beta^2/\tau^2) \rmd s \\
 &\le&  \int_{\frac12 \beta^2/\tau^2}^\infty e^{-\tau^2 (s-\frac12 \beta^2/\tau^2)} cG_y(x,s) \rmd s \\
 &\le&  c e^{\frac12 \beta^2} \int_0^\infty e^{-\tau^2 s} G_y(x,s) \rmd s \\
 &=&  c e^{\frac12 \beta^2} p_\tau(x;y) .
\end{eqnarray*}
 We then obtain\refq{est.Harnack-ellip}. \qed

\medskip
Let us remember that $y(\cdot)$ and $\Sigma$ were defined in Section \ref{ss.special}
and $P_\tau$ by\refq{def.PM}.

\begin{lem} (Caccioppoli's inequality for $P_\tau$).
\label{lem.Caccioppoli-PM}
 Let $P_\tau$ be defined by\refq{def.PM}. Let $\beta>0$.
 Then there exists $c>0$ such that for all $\tau >0$,
 if $B(\beta/\tau)\cap B(y(t);\frac1{\tau})=\emptyset$ we then have
\begin{equation}
\label{est.Caccioppoli-PM}
 \frac1{c}\int_{B(\frac{\beta}{4\tau})} \tau^2 P^2_\tau(x,t)  \rmd x \le 
 \int_{B(\frac{\beta}{2\tau})}  |\nabla P_\tau|^2(x,t) \rmd x \le
   c \int_{B(\frac{\beta}{\tau})} \tau^2 P^2_\tau(x,t)  \rmd x.
\end{equation}
\end{lem}
 Proof in Appendix.

\subsection{Comparison between $u_\tau$, $P_\tau$ and $p_\tau$}
\label{sec.comparison}
\begin{lem} (Comparison between $P_\tau$ and $p_\tau$).
\label{lem.CompPMptau}
 There exists $c>0$  such that for all $\tau >0$, $t\in \R$,
 if $x\not\in B(y(t); \frac2{\kappa^5 \tau})$ we then have
\begin{equation}
\label{comp.PMptau}
 \frac1{c} \tau^{3} P_\tau(x,t) \le  p_\tau(x,y(t)) \le  c\tau^{3} P_\tau(x,t) ,
\end{equation}
 where $\kappa$ is the constant of\refq{est.Gy} or\refq{est.ptau}.
\end{lem}
 Proof in Appendix.

\begin{lem} (Comparison between $u_\tau$ and $p_\tau$).
\label{lem.Comp-utau}
 Let $M>0$ and assume that $|\dot y|_\infty\le M$.
 Then there exist positive constants $C(M)$, $C_1(M)$, $\tau_0(M)$ such that
 for $\tau\ge\tau_0$, $t\in [0,T]$, $x\in\Om\sauf B(y(t),C_1/\tau)$,
 we have:
\begin{equation}
\label{est.utau}
 e^{-\tau^2(T+t)} u_\tau(x,t) \le  Ce^{-\tau\mu |t-\teta|} \tau^{-3} p_\tau(x,y(t)).
\end{equation}
\end{lem}
 Proof in Appendix.
\begin{lem}
\label{lem.dtutau}
 Let $t\in[0,T]$.
 Let $M>0$ and assume that $|\dot y|_\infty\le M$.
 Then there exist positive constants $C(M)$, $C_1(M)$, $\tau_0(M)$ such that 
 for $\tau\ge \tau_0$, $t\in [0,T]$ and $x\in\Om\sauf B(y(t),C_1/\tau)$ we have:
\begin{equation}
\label{est.dtutau}
 |\partial_t (e^{-\tau^2(t+T)} u_\tau(x,t))| \le C e^{-\tau\mu |t-\teta|} \tau^{-2} p_\tau(x,y(t)).
\end{equation}
\end{lem}
 Proof in Appendix.

\medskip
 Let us remember that $q_\tau$ is defined  by\refq{def.qtau}.
\begin{lem} (Estimate of $q_\tau$).
\label{lem.maj-qtau}
 Let $t\in[0,T]$.
 Let $M>0$ and assume that $|\dot y|_\infty\le M$.
 Then there exist positive constants $C(M)$, $C_1(M)$, $\tau_0(M)$ such that 
 for $\tau\ge \tau_0$, $t\in [0,T]$ and $x\in\Om\sauf B(y(t),C_1/\tau)$, we have
\begin{equation}
\label{est.qtau}
 |q_\tau(x,t)|   \le   C \tau^{-3} e^{-\tau\mu |t-\teta|} |x-y(t)| \, p_\tau(x,y(t)) .
\end{equation}
\end{lem}
 Proof in Appendix.

\subsection{Estimates of special function in $D_t$}
\begin{lem} (Estimates of $P_\tau$ in $D_t$).
\label{lem.nablaPtau-Dt} 
 Let $t\in[0,T]$.
 Then, there exists $c\ge 1$ such that for all $\tau> \frac{12}{\kappa^5 \dtyt}$, we have
\begin{equation}
\label{est.ptau-Dt}
 \frac1{c}\int_{D_t}  \tau^{-4} p^2_\tau(x,y(t))  \rmd x \le \int_{D_t}  |\nabla P_\tau(x,t)|^2 \rmd x
 \le c\int_{D_t}  \tau^{-4} p^2_\tau(x,y(t))  \rmd x.
\end{equation}
\end{lem}
 Proof in Appendix.

\begin{lem} (Estimate of $\nabla q_\tau$ in $D_t$).
\label{lem.maj-dqtau}
 Let $M>0$ and assume that $|\dot y|_\infty\le M$.
 Then there exist two positive constants $C_M$ and $\tau_0=\tau_0(\Sigma)$
 such that if $\tau> \tau_0$, $t\in[0,T]$, then
\begin{equation}
\label{est.dqtau-Dt}
 \int_{D_t}  |\nabla q_\tau|^2(x,t) \rmd x \le
  C_M \tau^{-4} e^{-2\tau\mu |t-\teta|} \int_{D_t}  |x-y(t)|^2 |p_\tau(x,y(t))|^2 \rmd x.
\end{equation}
\end{lem}
 Proof in Appendix.

\begin{lem}
\label{lem.(x-y)ptau}
 There exist positive constant $C$,  $\tau_0(\Sigma)$ such that for $\tau>\tau_0$, $t\in [0,T]$,
 we have
\begin{equation}
\label{est.(x-y)ptau}
 \int_{D_t} |x-y(t)|^2 |p_\tau(x,y(t))|^2 \rmd x  \le C \dtyt^2 \int_{D_t}  |p_\tau(x,y(t))|^2 \rmd x.
\end{equation}
\end{lem}
 Proof in Appendix.
 
\medskip
 Now we are ready to prove Lemma \ref{lem.R9}.

\subsection{Proof of Lemma  \ref{lem.R9}}
 We obtain\refq{minIinfty01} and\refq{majIinfty01} from\refq{minIinfty0},\refq{majIinfty0},
 \refq{est2.dqtau-Dt} of Lemma \ref{lem2.maj-dqtau},\refq{est.ptau-Dt} of Lemma  \ref{lem.nablaPtau-Dt}.
\qed

\section*{\textcolor{black}{Conclusion}}
\textcolor{black}{
 So we have proven the injectivity of $D\mapsto \Lambda_{v_0,{\bf a}}$ by extending the Dynamical Probe Method.
 We already know that, in the case where ${\bf a}$ is scalar and the background ${\bf b}=1$,
 the DPM is effective in reconstructing the inclusion $D$ from the Dirichlet-to--Neumann mapping $\Lambda_{v_0,{\bf a}}$
 even when $D_t$ has no regularity according to the space variable.
 But in the more general case where ${\bf b}$ is not constant the behaviour of the special functions $u_\tau$, $u_\tau^*$
 as $\tau$ tends to infinity is not so obvious anymore, which technically requires us to prove that
 the product $\nabla u_\tau \nabla u_\tau^*$ is positive not punctually but in a weaker sense,
 and with additional conditions.
 In our work the main new constraint that allows the uniqueness proof to work is on the geometry of D:
 some kind of uniform lipschitzian regularity of $D_t$, $t\in [0,T]$.
 By looking carefully at the various technical elements of the multiple Lemmas we can hope to improve
 this condition a little, perhaps by replacing it by a geometric constraint of the Holder type with coefficient in $(\frac12,1)$.
 The question of reconstructing $D$ from $\Lambda_{v_0,{\bf a}}$ remains delicate for two reasons.
 First, the negative term $C_M\dtyt^2$ in\refq{minIinfty01} forces the curve to be partly sufficiently 
 close enough to the inclusion to obtain a good lower bound of the indicator function $I_\infty(\tau)$,
 which complicates a strategy for detecting the unknown $D$.
 Then, Runge's method allows only a theoretical reconstruction. 
 Nevertheless, the reconstruction of points of $D$ sufficiently close to the lateral boundary of the cylinder
 becomes possible, and this without the use of the Runge approximation. 
 However, such a study would burden the article.\\
 Another question is to be able to weaken the condition that $t\mapsto D_t$ is lipschitzian. It is open.
 }

\section*{Appendix}
 \ssl {Proof of Lemma \ref{est-lem.Iinfty}}.
 We put
\begin{eqnarray*}
  X_1 &:=& \int_{\Om\times [0,T]}  ({\bf b}^{-1} - {\bf a}^{-1}) ({\bf b}\nabla u_\tau)^2 
 \rmd x e^{-2\tau^2 (T+t)} \rmd t ,\\
  X_2 &:=& \int_{\Om\times [0,T]} ({\bf a}-{\bf b})\, (\nabla u_\tau)^2 \, \rmd x
   e^{-2\tau^2 (T+t)} \rmd t,
\end{eqnarray*}
\begin{eqnarray}
\label{rel.w=u-v}
 w_\tau &:=& v_\tau - u_\tau ,\\
\label{def.Psi}
 \Psi_\tau &:=& ({\bf a}-{\bf b})\nabla v_\tau + {\bf b} \nabla w_\tau 
 =  {\bf a}\nabla v_\tau - {\bf b} \nabla u_\tau\\
\nonumber  &=&  ({\bf a} - {\bf b}) \nabla u_\tau + {\bf a}\nabla w_\tau,
\end{eqnarray}
\begin{eqnarray*}
  B_1 &:=&  \int_{\Om\times [0,T]}  {\bf a}^{-1} (\Psi_\tau)^2 \, \rmd x \, e^{-2\tau^2 (T+t)} 
  \rmd t \\
  B_2 &:=&  \int_{\Om\times [0,T]} {\bf a} (\nabla w_\tau)^2 \,  \rmd x  \, e^{-2\tau^2 (T+t)}
  \rmd t,\\
  B_3 &:=&   \int_{\Om\times [0,T]} \tau^2 w_\tau^2  \, \rmd x \, e^{-2\tau^2 (T+t)} \rmd t,
\end{eqnarray*}
 and
\begin{eqnarray*}
  R_1 &:=& \int_\Om \left[w_\tau  u_\tau^*\right]^T_0 \rmd x,\\
  R_2 &:=& \int_{\Om\times [0,T]}  ({\bf a}-{\bf b})  \nabla v_\tau
  \nabla (e^{\tau^2 (T+t)} u_\tau^* - e^{-\tau^2 (T+t)} u_\tau) \,  \rmd x \, e^{-\tau^2 (t+T)} \rmd t  \\
  &=& \int_{\Om\times [0,T]}  ({\bf a}-{\bf b}) \nabla v_\tau \cdot 
  (\nabla q_\tau^*(x,t) - \nabla q_\tau(x,t)) \,  \rmd x \, e^{-\tau^2 (t+T)} \rmd t ,\\
  R_3 &:=& \frac12 \int_\Om \left[w_\tau^2  e^{-2\tau^2 (T+t)} \right]^T_0 \rmd x.
\end{eqnarray*}
%

\medskip
 {\bf Step 1.}  We prove that
\begin{eqnarray}
\label{minIinfty1}
  I_\infty(\tau)   &=& X_1 + B_1 + B_3  +  R_1 + R_2 + R_3 ,
\end{eqnarray}
\begin{eqnarray}
\label{majIinfty1}
  I_\infty(\tau)   &=& X_2 - B_2 - B_3 +  R_1 + R_2 - R_3  .
\end{eqnarray}
 \medskip
 From\refq{d.Iinfty} we have
\begin{equation}
\label{r1.Infty}
 I_\infty(\tau)   =  \int_{\Om\times [0,T]}({\bf a}-{\bf b}) \nabla v_\tau
  \nabla u_\tau \,  \rmd x  e^{-2\tau^2 (T+t)} \rmd t + R_1 + R_2.
\end{equation}
 1. We put
\begin{eqnarray*}
  A_1 &:=&  \int_{\Om\times [0,T]} {\bf a}^{-1} \Psi_\tau  \cdot ({\bf a}-{\bf b}) \nabla u_\tau \,  \rmd x
   \,  e^{-2\tau^2 (T+t)} \rmd t,\\
  A_2 &:=&  \int_{\Om\times [0,T]} \nabla w_\tau  \Psi_\tau \,  \rmd x \,  e^{-2\tau^2 (T+t)} \, \rmd t.
\end{eqnarray*}
 Then, since  $({\bf a} - {\bf b}) \nabla u_\tau  =\Psi_\tau   - {\bf a}\nabla w_\tau$,
 we then have $ A_1  =  B_1 - A_2$.

\medskip
 By integration by parts we have
\begin{eqnarray}
\nonumber
  A_2  &=&  - \int_{\Om\times [0,T]}  w_\tau  \, \dive \Psi_\tau \,  \rmd x  e^{-2\tau^2 (T+t)}  \rmd t \\
\label{v.A2}  &=&  - \int_{\Om\times [0,T]}  w_\tau  \, \partial_t  w_\tau\,  \rmd x  e^{-2\tau^2 (T+t)} \rmd t 
 =  - B_3 - R_3.
\end{eqnarray}
 We thus have
\begin{equation}
\label{r.A1A2}
 A_1  =  B_1 + B_3 + R_3.
\end{equation}
 For any $3\times 3$ real matrix ${\bf m}$ we have ${\bf m}\nabla u_\tau \cdot   \nabla u_\tau=
 {\bf m}_S \nabla u_\tau \cdot   \nabla u_\tau$. 
 Then, thanks to 
\begin{equation}
\label{r2.vtau-utau}
 \nabla v_\tau = {\bf a}^{-1} \Psi_\tau + {\bf a}^{-1}{\bf b} \nabla u_\tau,
\end{equation}
 we obtain\refq{minIinfty1} from\refq{r1.Infty} and\refq{r.A1A2}.

 2. We consider\refq{r1.Infty} again. Thanks to\refq{rel.w=u-v} then to\refq{def.Psi}
 we have
\begin{eqnarray*}
 ({\bf a}-{\bf b}) \nabla v_\tau  \nabla u_\tau  = 
 ({\bf a}-{\bf b}) \nabla u_\tau\nabla u_\tau + {\bf a}\nabla w_\tau \nabla u_\tau
  - {\bf b}\nabla w_\tau \nabla u_\tau \\
  = ({\bf a}-{\bf b}) \nabla u_\tau\nabla u_\tau + {\bf a}\nabla w_\tau (\nabla v_\tau
  - \nabla w_\tau)  - {\bf b}\nabla w_\tau \nabla u_\tau \\
  =  ({\bf a}-{\bf b}) \nabla u_\tau\nabla u_\tau - {\bf a} \nabla w_\tau\nabla w_\tau
   +  \nabla w_\tau \Psi_\tau.
\end{eqnarray*}
 Hence
\begin{eqnarray*}
 I_\infty(\tau) &=&  \int_{\Om\times [0,T]}({\bf a}-{\bf b}) (\nabla u_\tau)^2 \,  \rmd x
  e^{-2\tau^2 (T+t)} \rmd t - B_2 + A_2 +  R_1 + R_2,
\end{eqnarray*}
 which gives\refq{majIinfty1} with the help of\refq{v.A2}.
\qed

\medskip
 {\bf Step 2.} 
 We put
\begin{equation}
\label{def.X0}
  X_0 := \int_{D}  e^{-2\tau\mu |t-\teta|} |\nabla P_\tau(x,t)|^2 \, \rmd x \rmd t
\end{equation}
 and
\begin{eqnarray*}
 R_4 &:=& \frac12 \int_\Om e^{4\tau^2 T} |u_\tau^*(T)|^2   \rmd x  + \frac12 \int_\Om e^{2\tau^2 T}
   |u_\tau^*(0)|^2  \rmd x \\
  &&  + 2\int_\Om  e^{-2\tau^2 T}  |u_\tau(0)|^2 \rmd x +  2\int_\Om  e^{-2\tau^2 T} |v_0|^2 \rmd x,\\
 R_5 &:=& \int_{D}  |\nabla q_\tau|^2   \,  \rmd x \rmd t, \quad
 R_5^* := \int_{D}  |\nabla q_\tau^*|^2   \,  \rmd x \rmd t.
\end{eqnarray*}
Thanks to\refq{h.gamma} and to Assumption (H0b) we have the following estimates:
\begin{eqnarray}
\label{minIinfty2}
  I_\infty(\tau)   &\ge&  C  X_0 + \frac12 B_1 + B_3 - 2R_4 - \frac1C (R_5+R_5^*) ,
\end{eqnarray}
\begin{eqnarray}
\label{majIinfty2}
 I_\infty(\tau)   &\le&  \frac1C X_0 - \frac12 B_2 - B_3 + 2R_4+  \frac1C (R_5+R_5^*) ,
\end{eqnarray}
 for some $C\in (0,1)$.
%
%

 Proof. Thanks to Cauchy-Minkovski inequality and to the definition\refq{def.wtau} we have
\begin{eqnarray}
\nonumber R_1+R_3 &=& 
 \int_\Om (w_\tau  u_\tau^* + \frac12 w_\tau^2  e^{-4\tau^2 T})|_{t=T}   \rmd x \\
\nonumber  && -  \int_\Om (w_\tau  u_\tau^* + \frac12 w_\tau^2  e^{-2\tau^2 T})|_{t=0}  \rmd x   \\
\label{est1.R13}  &\ge& 
  -R_4. 
\end{eqnarray}
 Similarly we have
\begin{eqnarray}
 \label{est2.R13} 
  R_1- R_3 &\le& R_4. 
\end{eqnarray}
 We observe that, thanks to\refq{minIinfty1} and\refq{majIinfty1},
\begin{eqnarray}
\label{val.X1}
 X_1 = I_\infty(\tau) - B_1 - B_3 - R_1 - R_2 - R_3,\\
\label{val.X2}
 X_2 =  I_\infty(\tau) + B_2 + B_3 - R_1 - R_2 + R_3.
\end{eqnarray}
 Thanks to\refq{r2.vtau-utau} again we have
\begin{eqnarray}
\nonumber
 |R_2| &\le & \int_{\Om\times [0,T]} e^{-\tau^2 (T+t)} |{\bf a}-{\bf b}| |{\bf a}^{-1}| |\Psi_\tau|\,
   |\nabla q_\tau^* - \nabla q_\tau|  \, \rmd x \rmd t  \\
\nonumber
 && + \int_{\Om\times [0,T]} e^{-\tau^2 (T+t)} |{\bf a}-{\bf b}|  |{\bf a}^{-1}| \, |{\bf b}|
  |\nabla u_\tau|\,   |\nabla q_\tau^* - \nabla q_\tau| \,  \rmd x \rmd t  \\
 \label{est.R2}
  &\le &  \frac12 B_1  +  \frac12 X_1 + C(R_5+R_5^*).
\end{eqnarray}
 From\refq{val.X1} and\refq{est.R2} we get
\begin{equation}
\label{est2.R2}
 |R_2| \le   \frac12 I_\infty(\tau) -  \frac12(B_3+R_3+R_1+R_2) + C (R_5+R_5^*) .
\end{equation}
 Estimates\refq{minIinfty1},\refq{est1.R13} and\refq{est2.R2} imply
\begin{eqnarray}
\label{minIinfty02}
 I_\infty(\tau)   &\ge&  \frac12 X_1 + \frac12 B_1 + B_3 - R_4 - C (R_5+R_5^*) .
\end{eqnarray}
 By using\refq{def.qtau},\refq{def.qtau*}, (H0b), and the basic estimate
 $a^2\ge \frac12 (a+b)^2 - b^2$, we have
\begin{eqnarray*}
 X_1 & \ge & \int_{D} \delta_1 e^{-2\tau\mu |t-\teta|} |\nabla P_\tau(x,t)|^2 \, \rmd x  \rmd t \\
  && - \int_{\Om\times [0,T]} |{\bf b}| \, |{\bf a}^{-1}| |{\bf a}-{\bf b}|\,
   |\nabla q_\tau|^2 \, \rmd x  \rmd t  \\
  &\ge&  C X_0 - \frac1C  R_5,
\end{eqnarray*}
 for some $C\in (0,1)$. Then with\refq{minIinfty02} we obtain\refq{minIinfty2}.

 Similarly, by using\refq{majIinfty1},\refq{def.X0},\refq{val.X2} we obtain\refq{majIinfty2}.

\medskip
 {\bf Step 3.} 
  We prove that for $\tau>\mu+1$ we have
\begin{eqnarray}
\label{maj.R4}
 |R_4| &\le& (2\|v_0\|^2_{L^2(\Om)} +5d_\Om) e^{-\tau\mu \min(T-\teta,\teta)}.
\end{eqnarray}
 Proof.
 Firstly, we have
\begin{eqnarray*}
 0 \le u_\tau(x,0) &=& \int_0^\infty \int_{\R^3} e^{\tau^2 (T-s)} e^{-\tau\mu |\teta+s|} m(y,-s)
  G_{y}(x,s)  \rmd y \rmd s\\
 &\le & e^{\tau^2 T}  e^{-\tau\mu \teta} \int_0^\infty \int_{\R^3}  e^{-\tau^2 s} G_{y}(x,s)  \rmd y \rmd s
  = \frac 1{\tau^2} e^{\tau^2 T}  e^{-\tau\mu \teta} .
\end{eqnarray*}
 Here we used the notorious relation
\begin{equation} 
\label{val.intGy}
\int_{\R^3} G_{y}(x,s)  \rmd y =1.
\end{equation} 
 Hence
\begin{equation}
\label{est.u0}
   0\le e^{-\tau^2 T} u_\tau(x,0) \le  e^{-\tau\mu \teta}, \quad \tau\ge 1.
\end{equation}
 Similarly we have
\begin{equation}
\label{est.u*T}
   0\le e^{2\tau^2 T} u_\tau^*(x,T) \le  e^{-\tau\mu (T-\teta)} , \quad \tau\ge 1.
\end{equation}
 Secondly, since $\tau>\mu+1>1$ we have
\begin{eqnarray}
\nonumber
 0 \le u^*_\tau(x,0) &=& \int_0^\infty \int_{\R^3} e^{-\tau^2 (T+s)} e^{-\tau\mu |\teta-s|}
  m(y,s) G_{y}(x,s)  \rmd y \rmd s\\
\nonumber
 &\le & e^{-\tau^2 T} e^{-\tau\mu \teta} \int_0^\infty  e^{-(\tau^2-\tau\mu) s}  \int_{\R^3} G_{y}(x,s)  \rmd y \rmd s\\
 \label{est.u*0}
 &=& \frac 1{\tau^2-\tau\mu} e^{-\tau^2 T} e^{-\tau\mu \teta} \le e^{-\tau^2 T} e^{-\tau\mu \teta}.
\end{eqnarray}
 From\refq{est.u*0},\refq{est.u*T},\refq{est.u0}, we obtain for $\tau>\mu+1>1$:
$$
 R_4 \le 2\|v_0\|^2_{L^2(\Om)}  e^{-2\tau^2 T} + 2d_\Om e^{-\tau\mu (T-\teta)}  + 3d_\Om e^{-\tau\mu \teta},
$$
 which implies\refq{maj.R4}.
 
 Estimates\refq{minIinfty0} and\refq{majIinfty0} come immediately from\refq{minIinfty2},\refq{majIinfty2},\refq{maj.R4}
 and the fact that $B_j\ge 0$ for $j=1,2,3$.
\qed

\medskip
\ssl{Proof  of Lemma \ref{lem.Caccioppoli-PM}}.
 We observe that for all $t$, the function $P_\tau(\cdot;t)$ is the unique solution in $H^1(\R^3)$
 of
\begin{equation}
\label{eq.PM}
 (- \dive\, ( {\bf b} \nabla \cdot) + \tau^2)P_\tau(\cdot;t) = m_\tau(\cdot,t).
\end{equation}
 Let $\phi\in C^1(\R;[0,1])$ with $\phi(r)=1$ for $|r|\le 1/2$ and $\phi(r)=0$ for $|r|\ge 1$.
 Put $\psi(x)=\phi(\tau(x-x_0)/\beta)$ where $x_0$ is the center of the ball $B(\beta/\tau)$.
 We multiply\refq{eq.PM} by $P_\tau(\cdot,t)\psi^2$ and integrate it over $\Om$.
 Since $\supp(\psi) \cap \supp(m_\tau(\cdot,t))$ has Lebesgue measure zero, we then have
\begin{equation}
\label{r.ptauphi}
 \int_\Om  [{\bf b} (\nabla P_\tau(\cdot,t))^2 \psi^2
  + 2{\bf b}\nabla P_\tau(\cdot,t)  \psi \;  P_\tau(\cdot,t) \nabla \psi
   + \tau^2 P_\tau^2(\cdot,t) \psi^2 ] =0.
\end{equation}
 Then, from Cauchy-Minkovski's inequality,
\begin{eqnarray*}
 \int_\Om  [{\bf b} (\nabla P_\tau(,t))^2\psi^2 + \tau^2 P_\tau^2 (\cdot,t) \psi^2]
  \le \int_\Om  |2{\bf b} \nabla P_\tau(\cdot,t) \psi \; P_\tau(\cdot,t) \nabla \psi | \\
 \le   \int_\Om [\frac12 {\bf b} (\nabla P_\tau(\cdot,t))^2 \psi^2
  + 2{\bf b} P_\tau^2(\cdot,t) (\nabla \psi)^2].
\end{eqnarray*}
 Thus, for some $C'>0$, 
$$
 \int_\Om  [|\nabla P_\tau(,t)|^2 + \tau^2  P_\tau^2(\cdot,t)] \psi^2(x) \rmd x \le
   C' \int_\Om P_\tau^2(\cdot,t) |\nabla \psi|^2(x) \rmd x.
$$
 Since $\supp \psi\subset B(\beta/\tau)$ with $|\nabla\psi(x)|\le \frac{\tau}{\beta} \max|\phi'|$,
 $\psi\ge 0$, and $\psi= 1$ in $B(\frac{\beta}{2\tau})$, we then have
$$
 \int_{B(\frac{\beta}{2\tau})}  |\nabla P_\tau(\cdot,t)|^2(x) \rmd x \le
   C''\tau^2\int_{B(\frac{\beta}{\tau})} P_\tau^2(\cdot,t) \rmd x,
$$
 which proves the second inequality in\refq{est.Caccioppoli-PM}.\\
 From\refq{r.ptauphi} and thanks to Cauchy-Minkovski's inequality we have also
\begin{eqnarray*}
 \int_\Om  [{\bf b} (\nabla P_\tau(\cdot,t))^2 \psi^2 + \tau^2  P_\tau^2(\cdot,t)  \psi^2] \le
 \int_\Om  |2{\bf b}\nabla P_\tau(\cdot,t)  \nabla\psi \; P_\tau(\cdot,t)  \psi | \\
 \le  \int_\Om [\frac2{\tau^2} \gamma_\infty^2 |\nabla P_\tau(\cdot,t) |^2 |\nabla\psi|^2
  + \frac12 \tau^2 P_\tau^2(\cdot,t) \psi^2].
\end{eqnarray*}
 Thus, 
$$
 \int_\Om  \tau^2  P_\tau^2(\cdot,t) \psi^2(x) \rmd x \le
   C \tau^{-2} \int_\Om |\nabla P_\tau(\cdot,t) |^2 |\nabla \psi|^2(x) \rmd x.
$$
 We then obtain
$$
 \int_{B(\frac{\beta}{2\tau})}  \tau^2  P_\tau^2(\cdot,t)   \rmd x \le
   C' \int_{B(\frac{\beta}{\tau})} |\nabla P_\tau(\cdot,t) |^2 \rmd x
$$
 which proves the first inequality in\refq{est.Caccioppoli-PM} with $\beta$ replaced by $2\beta$.
 \qed

\medskip
 \ssl {Proof of Lemma \ref{lem.CompPMptau}}.
 Since $G_x(y,s) =G_y(x,s)$ and thanks to\refq{est.Harnack-parab} with $r=1/\tau$,
 we have for all $x\not\in B(y(t),2/\tau)$:
\begin{eqnarray*}
 \int_{B(y(t),1/\tau)} G_y(x,s)  \rmd y &=& \int_{B(y(t),1/\tau)} G_x(y,s)  \rmd y
 \le  |B(1/\tau)| \max_{B(y(t),1/\tau)} G_x(\cdot,s) \\
  &\le& c \tau^{-3} G_x(y(t),s+\frac1{2\tau^2}) = c \tau^{-3} G_{y(t)}(x,s+\frac1{2\tau^2}).
\end{eqnarray*}
 Then, since $\tau|x-y(t)|\ge 2/\kappa^5\ge 2$, since $m_\tau\le 1$ and supp $m_\tau=B(y(t),\frac1{\tau})$
 we have
\begin{eqnarray}
\label{est-proof.Ptau}
 P_\tau(x,t) &\le &    \int_0^{\infty} e^{-\tau^2 s}  \int_{B(y(t),\frac1{\tau})} G_y(x,s)  \rmd y \rmd s \\
\nonumber
  &\le &  c \tau^{-3} \int_0^{\infty} e^{-\tau^2 s} G_{y(t)}(x,s+\frac1{2\tau^2}) \rmd s \\
\nonumber
  & = &  c \tau^{-3} \int_{\frac1{2\tau^2}}^{\infty} e^{-\tau^2 (s-\frac1{2\tau^2})} G_{y(t)}(x,s) \rmd s \\
\nonumber
  &\le & c'  \tau^{-3} \int_0^{\infty}  e^{-\tau^2 s}  G_{y(t)}(x,s)  \rmd s  = c'  \tau^{-3} \; p_\tau(x;y(t)).
\end{eqnarray}
 We obtain the first inequality of\refq{comp.PMptau}. Let us prove the second one.
 Since $m_\tau\ge 1/2$ in $B(y(t),\frac1{2\tau})$ we then have
\begin{eqnarray*}
 P_\tau(x,t) &\ge &  \frac12 \int_0^{\infty} e^{-\tau^2 s}  \int_{B(y(t),\frac1{2\tau})} G_y(x,s)  \rmd y \rmd s \\
  &\ge &  c \tau^{-3}  \int_0^{\infty} e^{-\tau^2 s}  \inf_{y\in B(y(t),\frac1{2\tau})} G_y(x,s) \rmd s.
\end{eqnarray*}
 By applying\refq{est.Harnack-parab} with $r=1/\tau$ and observing that $G_y(x,s)=G_x(y,s)$
 we then have for all $x\not\in B(y(t),2/\tau)$:
\begin{eqnarray}
\nonumber
 P_\tau(x,t)  &\ge &  c \tau^{-3}  \int_0^{\infty} e^{-\tau^2 s}  \inf_{y\in B(y(t),\frac1{2\tau})} G_y(x,s) \rmd s \\
\nonumber
  &\ge &  c \tau^{-3}  \int_0^{\infty} e^{-\tau^2 s} G_{y(t)}(x,s-\frac1{2\tau^2}) \rmd s \\
\nonumber
  &=&  c \tau^{-3} \int_{\frac1{2\tau^2}}^{\infty} e^{-\tau^2 (s+\frac1{2\tau^2})} G_{y(t)}(x,s) \rmd s \\
\label{min0.Ptau}
  &= & c'  \tau^{-3} \left( p_\tau(x,y(t)) - \int_0^{\frac1{2\tau^2}}  e^{-\tau^2 s}  G_{y(t)}(x,s)  \rmd s\right),
\end{eqnarray}
 where $c'>0$.
 We put $R := \int_0^{\frac1{2\tau^2}}  e^{-\tau^2 s}  G_{y(t)}(x,s)  \rmd s$. 
 Thanks to\refq{est.Gy} and\refq{est.ptau} we have
\begin{eqnarray*}
 R &\le& \int_0^{\frac1{2\tau^2}}  \frac{e^{- \frac{\kappa^2|x-y(t)|^2}{4s}}}{\kappa s^{3/2}}  \rmd s
  \le  \sqrt2 \kappa^{-1}\tau \int_1^\infty e^{- \kappa^2|x-y(t)|^2 \tau^2 r/2}  \rmd r \\
 &=&  2\sqrt2 \kappa^{-3}\tau^{-1} |x-y(t)|^{-2} e^{- \kappa^2|x-y(t)|^2 \tau^2/2} \\
 &\le& \frac12 p_\tau(x,y(t)) \: \frac2{\kappa^5 \tau|x-y(t)|} \exp(-\kappa^{-1}\tau|x-y(t)|(\kappa^2\tau|x-y(t)|/2-1)).
\end{eqnarray*}
 Since $\tau|x-y(t)| \ge 2/\kappa^5>2/\kappa^2$, we then have $R\le \frac12  p_\tau(x,y(t))$.
 Hence
\begin{eqnarray*}
 P_\tau(x,t)  &\ge &  \frac12 c' \tau^{-3} p_\tau(x,y(t)).
\end{eqnarray*}
 The conclusion follows.
 \qed

\medskip
 \ssl {Proof of Lemma \ref{lem.Comp-utau}}.
 Let us observe that
\begin{equation}
\label{est.expteta}
 e^{-\tau\mu |t-\teta-s|} \le e^{-\tau\mu |t-\teta|} e^{\tau\mu s} \quad s>0,\; t\in\R .
\end{equation}
 Hence
\begin{eqnarray} 
\nonumber
  e^{-\tau^2(T+t)} u_\tau(x,t) &=& \int_{0}^\infty  e^{-\tau^2 s} e^{-\tau\mu |t-\teta-s|}
   \int_{\R^3} m_\tau(y,t-s) G_{y}(x,s)  \rmd y \rmd s \\
\label{est.utau-H}
  &\le &  e^{-\tau\mu |t-\teta|}  H,
\end{eqnarray}
 where we put
\begin{eqnarray} 
\label{def.H}
 H &:=& \int_{0}^\infty  e^{-(\tau^2 -\tau\mu) s}  \int_{B(y(t-s),1/\tau)} G_{y}(x,s)  \rmd y \rmd s
 \equiv H_1+H_2
\end{eqnarray}
 with
\begin{eqnarray*} 
 H_1 &:=& \int_{s>\lambda/\tau}  e^{-(\tau^2 -\tau\mu) s}  \int_{B(y(t-s),1/\tau)} G_{y}(x,s)  \rmd y \rmd s,\\
 H_2 &:=& \int_{0}^{\lambda/\tau}  e^{-(\tau^2 -\tau\mu) s}  \int_{B(y(t-s),1/\tau)} G_{y}(x,s)  \rmd y \rmd s,
\end{eqnarray*}
 and where $\lambda := 2|x-y(t)|/\kappa$. 
 We put also $M':=M\lambda+1$, $C_1=\max(1,8\kappa^{-7}(1+M^2 d_\Om^2))$.
 Since $|y(t-s)-y(t)|\le Ms$, we then have $B(y(t-s),1/\tau)\subset B(y(t),Ms+1/\tau)$
 and so
\begin{eqnarray*} 
 H_2 &\le &  e^{\mu\lambda}  \int_{0}^{\lambda/\tau}  e^{-\tau^2 s}  \int_{B(y(t),Ms+1/\tau)}
  G_{y}(x,s)  \rmd y \rmd s \\
 &\le &  e^{\mu\lambda}  \int_{0}^{\lambda/\tau}  e^{-\tau^2 s}  \int_{B(y(t),M'/\tau)}
  G_{y}(x,s)  \rmd y \rmd s.
\end{eqnarray*}
 Since $\tau\ge 2\kappa^{-1}M$, $|x-y(t)|\ge 1/\tau$, we then have $|x-y(t)|\ge M'/\tau$
 and so we can apply\refq{est.Harnack-parab} where $x$ and $y$ are exchanged and
 with $r=M'/(2\tau)$. Hence
\begin{eqnarray*} 
 H_2 &\le &  ce^{\mu\lambda}  \int_{0}^\infty  e^{-\tau^2 s} |B(y(t),M'/\tau)| G_{y(t)}(x,s+{M'}^2/(2\tau^2)) \rmd s\\
 &\le & c e^{\mu\lambda}  \int_{0}^\infty  e^{-\tau^2 s}  (2M')^3 |B(y(t),1/(2\tau))|
  G_{y(t)}(x,s+{M'}^2/(2\tau^2)) \rmd s\\
 & = & c e^{\mu\lambda}   (2M')^3 \int_{{M'}^2/(2\tau^2)}^\infty  e^{-\tau^2 (s-{M'}^2/(2\tau^2))}
   |B(y(t),1/(2\tau))|  G_{y(t)}(x,s) \rmd s\\
 &= & c e^{\mu\lambda+{M'}^2/2}  {M'}^3 \tau^{-3} (p_\tau(x,y(t)) - R),
\end{eqnarray*}
 with $R := \int_{0}^{{M'}^2/(2\tau^2)} e^{-\tau^2 s}  G_{y(t)}(x,s) \rmd s$ and
 $c$ is the constant\refq{est.Harnack-parab}.
 Thanks to\refq{est.Gy} we have
\begin{eqnarray*} 
 R &\le&  \int_0^{\frac{{M'}^2}{2\tau^2}}  \frac{e^{-\frac{\kappa^2|x-y(t)|^2}{4s}}}{\kappa s^{3/2}}  \rmd s
  \le  \sqrt2 \kappa^{-1}\tau {M'}^{-1} \int_1^\infty e^{-\kappa^2|x-y(t)|^2 \tau^2 r/(2{M'}^2)}  \rmd r \\
  &= & 2\sqrt2 {M'} \kappa^{-3}\tau^{-1} |x-y(t)|^{-2} e^{-\kappa^2|x-y(t)|^2 \tau^2 /(2{M'}^2)}  \\
  &\le& \frac12 p_\tau(x,y(t)) \: \frac{M'}{\kappa^5 \tau|x-y(t)|}
  \exp(-\kappa^{-1}\tau|x-y(t)|({M'}^{-2}\kappa^2\tau|x-y(t)|/2-1)).
\end{eqnarray*}
 Since $\tau|x-y(t)|\ge C_1$ then ${M'}^{-2}\kappa^{2}\tau|x-y(t)|\ge 2$ and
 ${M'}^{-1}\kappa^5 \tau|x-y(t)|\ge 1$. Hence $R\le\frac12 p_\tau(x,y(t))$ and
\begin{eqnarray} 
\nonumber
 H_2 &\le & C {M'}^3 e^{2\kappa^{-1}\mu|x-y(t)|} \tau^{-3} p_\tau(x,y(t))\\
\nonumber
  &\le &C {M'}^3 e^{2\kappa^{-1}\mu d_\Om} \tau^{-3} p_\tau(x,y(t)) , \quad (x,t)\in\Om_{0,T}\\
\label{est.H2}
  &=& C'(M) \; \tau^{-3} p_\tau(x,y(t)), \quad (x,t)\in\Om_{0,T}.
\end{eqnarray}
 Let us estimate $H_1$. 
 Since $G_y(x,s)\le \kappa^{-1}s^{-3/2}$ and $\tau\ge 2\mu$ we then have
\begin{eqnarray*} 
 H_1 &\le & \kappa^{-1} \tau^{-3} \int_{s>\lambda/\tau}  s^{-3/2} e^{-\tau^2 s/2}  \rmd s\\
&\le & \kappa^{-1}   \tau^{-3}
 \left\{\begin{array}{ccc}
  \int_{s>\lambda/\tau}   (\lambda/\tau)^{-3/2} e^{-\tau^2 s/2}  \rmd s &=& 2 (\lambda/\tau)^{-3/2}\tau^{-2} e^{-\tau\lambda/2}\\
  \int_{s>\lambda/\tau}   s^{-3/2} e^{-\tau\lambda/2}  \rmd s &=& 2(\lambda/\tau)^{-1/2} e^{-\tau\lambda/2}\\
\end{array} \right. .
\end{eqnarray*}
 Hence
 $$ H_1 \le  2 \kappa^{-1}  \lambda^{-1}  \tau^{-3} e^{-\tau\lambda/2}.$$
 Thanks to\refq{est.ptau}  we then obtain
\begin{eqnarray} 
\label{est.H1}
 H_1 &\le &  \kappa^{-2} \tau^{-3} p_\tau(x,y(t)), \quad (x,t)\in\Om_{0,T}.
\end{eqnarray}
 Then, thanks to Lemma \ref{lem.CompPMptau} and from\refq{est.H1},\refq{est.H2},\refq{est.utau-H},
 the conclusion follows.
\qed

\medskip
 \ssl {Proof of Lemma \ref{lem.dtutau}}.
 We have
\begin{eqnarray*} 
  \partial_t (e^{-\tau^2(t+T)} u_\tau(x,t)) &=& Y_1 + Y_2,
\end{eqnarray*}
 with
\begin{eqnarray*} 
  Y_1 &:=& -\tau\mu \int_{s=0}^t  e^{-\tau^2 s} \sgn(t-\teta-s) e^{-\tau\mu |t-\teta-s|} \int_{\R^3} m_\tau(y,t-s)
  G_{y}(x,s)  \rmd y \rmd s, \\
  Y_2 &:=& \int_{s=0}^t  e^{-\tau^2 s}  e^{-\tau\mu |t-\teta-s|} \int_{\R^3} \partial_t m_\tau(y,t-s)
  G_{y}(x,s)  \rmd y \rmd s.
\end{eqnarray*}
 Let us estimate $Y_1$. We have
\begin{eqnarray*}
  |Y_1(x,t)| &\le & \tau\mu \int_{s=0}^\infty  e^{-\tau^2 s} e^{-\tau\mu |t-\teta-s|} \int_{\R^3}
   m_\tau(y,t-s) G_{y}(x,s)  \rmd y \rmd s \\
  &=& \tau\mu\;  e^{-\tau^2(t+T)} u_\tau(x,t).
\end{eqnarray*}
  Thanks to Lemma \ref{lem.Comp-utau} we obtain 
\begin{equation}
\label{est.Y1}
  |Y_1(x,t)| \le  C e^{-\tau\mu |t-\teta|} \tau^{-2} p_\tau(x,y(t)).
\end{equation}
 Let us estimate $Y_2$. Remember that $\supp m_\tau(\cdot,t) \subset \ove{B(y(t),1/\tau)}$
 and that
$$ |\partial_t m_\tau(y,t)| = \tau |\dot y(t) \nabla M_0(\tau(y-y(t)))| \le M\tau . $$
 Hence we have , as in the estimates of\refq{est-proof.Ptau} we obtain 
\begin{eqnarray*}
  |Y_2| &\le & C\tau  \int_{s=0}^\infty  e^{-\tau^2 s}  e^{-\tau\mu |t-\teta-s|} \int_{B(y(t-s),1/\tau)} G_{y}(x,s)
    \rmd y \rmd s\\
 &\le&   e^{-\tau\mu |t-\teta|} H,
\end{eqnarray*}
 where $H$ is defined by\refq{def.H}.
 Hence
\begin{equation}
\label{est.Y2} 
  |Y_2| \le  C' e^{-\tau\mu |t-\teta|} \tau^{-2}  p_\tau(x,y(t)).
\end{equation}
 From\refq{est.Y1},\refq{est.Y2} we obtain\refq{est.dtutau}.
\qed

\medskip
 \ssl {Proof of Lemma \ref{lem.maj-qtau}}.
 We write $q_\tau(x,t) = \int_0^{\infty}  e^{-\tau^2 s}  \int_{\R^3} (A-B) G_y(x,s)  \rmd y \rmd s$
 with
\begin{eqnarray*}
 A &\equiv& e^{-\tau\mu |t-\teta-s|} m_\tau(y,t-s), \\
 B &\equiv& e^{-\tau\mu |t-\teta|} m_\tau(y,t).
\end{eqnarray*}
 Let us observe that, since $e^{\tau\mu s}-1\le  \mu \tau s e^{\tau\mu s}$ and
 thanks to\refq{est.expteta}, then
\begin{eqnarray*}
 |A-B| &\le&  e^{-\tau\mu |t-\teta|} \Big( \mu \tau s e^{\tau\mu s} 1_{B(y(t-s),1/\tau)} \\
  && + M\tau s \max(1_{B(y(t),1/\tau)},1_{B(y(t-s),1/\tau)})\Big).
\end{eqnarray*}
 Hence
\begin{equation}
\label{estqtau-R1R2}
 |q_\tau(x,t)| \le   \tau e^{-\tau\mu |t-\teta|} (\mu R_1+M R_2)
\end{equation}
 with
\begin{eqnarray}
\label{def.R1}
 R_1 & := &   \int_0^{\infty}  e^{-\tilde\tau^2 s} \int_{B(y(t-s),1/\tau)}  sG_y(x,s)  \rmd y \rmd s, \\
\label{def.R2}
 R_2 & := &  \int_0^{\infty} e^{-\tau^2 s} \int_{\tilde B} s G_y(x,s)  \rmd y \rmd s,
\end{eqnarray}
 where $\tilde\tau:=\sqrt{\tau^2-\tau\mu}$ and $\tilde B:=B(y(t-s),1/\tau)\cup B(y(t),1/\tau)$.\\
 Let us put again $\lambda=2\kappa^{-1}|x-y(t)|$.
 We write $R_2 = R_{21}+R_{22}$ with
\begin{eqnarray*}
 R_{21} &:=& \int_0^{\lambda/\tau}  e^{-\tau^2 s} \int_{\tilde B} sG_y(x,s)  \rmd y \rmd s, \\
 R_{22} &:=& \int_{\lambda/\tau}^{\infty}  e^{-\tau^2 s} \int_{\tilde B} sG_y(x,s)  \rmd y \rmd s.
\end{eqnarray*}
 As for the estimate of $H_1$ in the proof of Lemma \ref{lem.Comp-utau} we have
\begin{eqnarray*}
 |R_{22}| &\le& 
   2\kappa^{-1} \tau^{-3}  \int_{\lambda/\tau}^{\infty} s^{-1/2} e^{-\tau^2 s/2}  \rmd s\\
 &\le & 2\kappa^{-1}\tau^{-3+1/2} \lambda^{-1/2}  \int_{\lambda/\tau}^{\infty} e^{-\tau^2 s/2}  \rmd s
 = 4\kappa^{-1}\tau^{-5+1/2} \lambda^{-1/2} e^{-\tau \lambda/2}  \\
 & = &  2\sqrt2 \kappa^{-1/2} \tau^{-5+1/2} |x-y(t)|^{-1/2} e^{-\tau\kappa^{-1}|x-y(t)|} .
\end{eqnarray*}
 Thanks to\refq{est.ptau} and since $\tau|x-y(t)|\ge 1$ we then have
\begin{eqnarray}
\nonumber
 |R_{22}| &\le& \kappa^{-5/2}\tau^{-4} (\tau |x-y(t)|)^{-1/2} |x-y(t)| p_\tau(x,y(t))\\
\label{est.R22-qtau}
 &\le & C \tau^{-4} |x-y(t)| p_\tau(x,y(t)).
\end{eqnarray}
 For $s\le \lambda/\tau$ we have $\tilde B\subset B(y(t),M'/\tau)$ with $M':=\lambda M+1$.
 Hence, as for the estimate of $H_2$ in the proof of Lemma \ref{lem.Comp-utau} and
 since $\tau|x-y(t)|\ge 2{M'}^2 \kappa^{-5}$ we then have
\begin{eqnarray}
\nonumber
 |R_{21}| &\le & \lambda\tau^{-1}  \int_{0}^{\lambda/\tau}  e^{-\tau^2 s}  \int_{B(y(t),M'/\tau)} G_y(x,s)  \rmd y \rmd s \\
\nonumber
 &\le &  \lambda\tau^{-1}  C(M)  \tau^{-3} p_\tau(x,y(t)) \\
 \label{est.R21-qtau} 
 &\le &  C(M) \tau^{-4} |x-y(t)| p_\tau(x,y(t)).
\end{eqnarray}
 From\refq{est.R22-qtau} and\refq{est.R21-qtau}, we obtain that for $\tau|x-y(t)|\ge 2{M'}^2 \kappa^{-5}$
 we have
\begin{eqnarray}
\label{est.R2-qtau}
 |R_2| &\le & C(M) \tau^{-4} |x-y(t)|  p_\tau(x,y(t)).
\end{eqnarray}
 Now, we estimate $R_1$ as $R_2$ by splitting the integral in\refq{def.R2} with
 $s<\lambda/\tau$ or $s>\lambda/\tau$.
 We observe that $\tilde\tau = \tau\sqrt{1-\mu\tau^{-1}}\ge \frac1{\sqrt2}\tau$.
 Hence $R_1=R_{11}+R_{12}$ with, since $\tilde\tau|x-y(t)|\ge \sqrt{2} \kappa^{-5}$ and $\tau\ge2\mu$,
\begin{eqnarray}
\nonumber
 |R_{12}| &\le & \tau^{-3+1/2}\lambda^{-1/2}  \int_{\lambda/\tau}^{\infty}  e^{-\tilde\tau^2 s/2}  \rmd s\\
\nonumber
 &=& 2\kappa^{-1} \tilde\tau^{-2}\tau^{-3+1/2} \lambda^{-1/2} e^{-\tau\lambda+\lambda\mu}  \\
\label{est.R12-qtau}
 &\le &  C \tau^{-4} |x-y(t)| p_\tau(x,y(t)).
\end{eqnarray}
 Finally, since $\tau|x-y(t)|\ge 2\kappa^{-5}$, we have, as for the estimate of $R_{21}$,
\begin{eqnarray}
\nonumber
 |R_{11}| &\le & \lambda\tau^{-1} e^{\lambda\mu}  \int_0^{\lambda/\tau} e^{-\tau^2 s} 
  \int_{B(y(t),1/\tau)} G_y(x,s)  \rmd y \rmd s\\
\label{est.R11-qtau}
 &\le &  C(M)  \tau^{-4} |x-y(t)| p_\tau(x,y(t)).
\end{eqnarray}
 Thanks to\refq{est.R12-qtau} and\refq{est.R11-qtau}, we obtain
\begin{eqnarray}
\label{est.R1-qtau}
 |R_1| &\le & C(M) \tau^{-4} |x-y(t)| p_{\tau}(x,y(t)).
\end{eqnarray}
 Thanks to\refq{est.R1-qtau} and\refq{est.R2-qtau},\refq{estqtau-R1R2}, we obtain\refq{est.qtau}
 for $\tau\ge 2\mu$, $t\in [0,T]$ and $x\in\Om\sauf B(y(t),C_1/\tau)$.
\qed

\medskip
 \ssl {Proof of Lemma \ref{lem.nablaPtau-Dt}}.
 We consider a familly of balls $B_i(1/\tau)$, $i\in I$, as in Lemma \ref{lem.geometrie}.
 By using\refq{est.Harnack-ellip} (with $\beta=6$),\refq{est.Caccioppoli-PM} and\refq{comp.PMptau},
 and by observing that $B_i(6/\tau)\cap B(y(t),2\kappa^5/\tau)=\emptyset$
 for $\tau >12\kappa^{-5}\dtyt^{-1}$, we can write
\begin{eqnarray}
\nonumber
 \int_{D_t} |\nabla P_\tau|^2 (x,t) \rmd x &\le&
 \sum_{i\in I} \int_{B_i(3/\tau)} |\nabla P_\tau|^2 (x,t) \rmd x \\
\nonumber
 &\le& C_1 \sum_{i\in I} \int_{B_i(6/\tau)} \tau^2 P_\tau^2 (x,t) \rmd x \\
\nonumber
 &\le &C_2\sum_{i\in I}   \int_{B_i(6/\tau)} \tau^{-4} p_\tau^2 (x,y(t)) \rmd x \\
\nonumber
 &\le& C_3 \sum_{i\in I} |B_i(6/\tau)| \tau^{-4}  p_\tau^2 (x_i,y(t)) \\
\nonumber
 &\le & C_4 \sum_{i\in I} |B_i(1/\tau)| \tau^{-4} \min_{B_i(1/\tau)} p_\tau^2 (\cdot,y(t)) \\
\nonumber
 &\le & C_5 \sum_{i\in I} \int_{B_i(1/\tau)} \tau^{-4} p_\tau^2 (x,y(t)) \rmd x \\
\label{maj.ptau-onDt}
 &\le& C_5 \int_{D_t} \tau^{-4} p_\tau^2 (x,y(t)) \rmd x.
\end{eqnarray}
 Hence, the second inequality of\refq{est.ptau-Dt} is proved. The proof of the first one is similar.
\qed

\medskip
 \ssl {Proof of Lemma \ref{lem.maj-dqtau}}.
 We put $C'_1(M)=C_1+6$, $C_2=\max(C'_1,12\kappa^{-5})$ where $C_1(M)$ is
 the constant in Lemma \ref{lem.maj-qtau}.
 We consider again the balls $B(1/\tau)$, $B(3/\tau)$, defined in Lemma \ref{lem.geometrie}.
 Thus
\begin{eqnarray}
\label{def.J-maj}
 J&:=& \int_{D_t} |\nabla q_\tau (x,t)|^2\rmd x \le  \sum_i \int_{B_i(3/\tau)} |\nabla q_\tau (x,t)|^2\rmd x.
\end{eqnarray}
 Let us fix $i$ and denote $B(3/\tau)=B_i(3/\tau)$.
 We consider again the functions $\phi\in C^1(\R;[0,1])$ and $\psi(x)=\phi(\tau(x-x_0)/6)$ where $x_0$
 is the center of a ball $B(6/\tau)$, as in the proof of Lemma \ref{lem.Caccioppoli-PM} (with $\beta=6$).\\
 Thanks to Lemma \ref{lem.dtutau}, there exists a positive constant $C(M)$
 such that for $\tau\ge2\mu$, $t\in [0,T]$, $x\in \Om\sauf B(y(t);C_1/\tau)$, we have
\begin{eqnarray}
\nonumber
 \left| (-\dive\,{\bf b}\nabla + \tau^2) q_\tau \: (x,t)\right| &=& \left|  \partial_t (e^{-\tau^2(t+T)}u_\tau(x,t)) \right|\\
\label{eq.qtau}
  &\le& C \tau^{-2} e^{-\tau\mu|t-\teta|}  p_\tau(x,y(t)).
\end{eqnarray}
 We observe that
\begin{eqnarray*}
 x\in \supp(\psi) = \ove{B(6/\tau)} \Rightarrow |x-y(t)| \ge |x_0-y(t))|-6/\tau 
 \ge \\
 \dtyt-6/\tau >C'_1/\tau - 6/\tau = C_1/\tau.
\end{eqnarray*}
 Hence we can multiply\refq{eq.qtau} by $q_\tau(x,t)\psi^2(x)$ and integrate it over $\Om$. 
 This implies
\begin{eqnarray*}
 \int_\Om  \left( {\bf b} (\nabla q_\tau(\cdot,t))^2 \psi^2 + 2{\bf b}\nabla q_\tau(\cdot,t)  \psi \; 
  q_\tau(\cdot,t) \nabla \psi + \tau^2  q_\tau^2(\cdot,t) \psi^2 \right) \\
 \le  C\tau^{-2}   e^{-\tau\mu |t-\teta|} \int_\Om  |q_\tau(\cdot,t)|  p_\tau(\cdot,y(t)) \psi^2.
\end{eqnarray*}
 Then, from Cauchy-Minkovski's inequality, and as in the proof of Lemma \ref{lem.Caccioppoli-PM},
 we obtain
\begin{eqnarray*}
 \int_\Om  (|\nabla q_\tau(,t)|^2 + \tau^2  q_\tau^2(\cdot,t)) \psi^2 
 &\le&  C \int_\Om q_\tau^2(\cdot,t) (\nabla \psi)^2 +  C  e^{-\tau\mu |t-\teta|} \tau^{-2} \cdot \\
 && \hspace{-0.4cm} \big(\int_\Om  |q_\tau(\cdot,t)|^2\psi^2\big)^{1/2}
  \big(\int_\Om |p_\tau(\cdot,y(t))|^2 \psi^2\big)^{1/2}.
\end{eqnarray*}
 Since $\supp \psi= \ove{B(6/\tau)}$ with $|\nabla\psi(x)|\le \tau \max|\phi'|/6$, $\psi\ge 0$,
 and $\psi= 1$ in $B(\frac3{\tau})$, we then have
\begin{eqnarray*}
 \int_{B(\frac3{\tau})}  |\nabla q_\tau(\cdot,t)|^2  &\le&  C\tau^2 \int_{B(\frac6{\tau})} q^2_\tau(\cdot,t) +
  C e^{-\tau\mu |t-\teta|} \tau^{-2} \cdot \\
  && \cdot \big(\int_{B(\frac6{\tau})}  |q_\tau(\cdot,t)|^2\big)^{1/2}
  \big(\int_{B(\frac6{\tau})} |p_\tau(\cdot,y(t))|^2 \big)^{1/2} .
\end{eqnarray*}
 Thanks to Lemma \ref{lem.maj-qtau} and by using $\tau^{-1}\le C_1|x-y(t)|$ for $x\in B(\frac6{\tau})$,
 we then have 
\begin{eqnarray}
\nonumber
 \int_{B(\frac3{\tau})}  |\nabla q_\tau(\cdot,t)|^2  &\le&
  C\tau^{-4}  e^{-2\tau\mu |t-\teta|}  \int_{B(\frac6{\tau})}  |x-y(t)|^2 |p_\tau(x,y(t))|^2 \rmd x  +\\ 
\nonumber
 && C\tau^{-5}  e^{-2\tau\mu |t-\teta|}   \cdot  \big(\int_{B(\frac6{\tau})}
   |x-y(t)|^2 |p_\tau(x,y(t))|^2 \rmd x \big)^{1/2} \cdot\\
\nonumber
  && \cdot \big( \int_{B(\frac6{\tau})} |p_\tau(x,y(t))|^2 \rmd x \big)^{1/2}  \\
\label{est.dqtauBall}
   &\le& C' e^{-2\tau\mu |t-\teta|} \tau^{-4}  \int_{B(\frac6{\tau})}  |x-y(t)|^2 |p_\tau(x,y(t))|^2 \rmd x.
\end{eqnarray}
 By putting\refq{est.dqtauBall} in\refq{def.J-maj} we obtain
\begin{equation}
\label{est1.dqtau}
 J \le C' e^{-2\tau\mu |t-\teta|} \tau^{-4}  \sum_i  \int_{B(\frac6{\tau})}  |x-y(t)|^2 |p_\tau(x,y(t))|^2 \rmd x.
\end{equation}
 Finally, as in\refq{maj.ptau-onDt} we have 
$$
  \sum_i  \int_{B(\frac6{\tau})}  |x-y(t)|^2 |p_\tau(x,y(t))|^2 \rmd x
 \le C  \int_{D_t}  |x-y(t)|^2 |p_\tau(x,y(t))|^2 \rmd x.
$$
 This with\refq{est1.dqtau} prove\refq{est.dqtau-Dt}.
\qed

\medskip
\ssl {Proof of Lemma \ref{lem.(x-y)ptau}}.
 We can assume that $D_t\neq \emptyset$.
 We put $\lambda=2\kappa^{-2} \dtyt$  and
\begin{eqnarray*}
 J &:=& \int_{D_t} |x-y(t)|^2 |p_\tau(x,y(t))|^2 \rmd x = J_1+J_2,\\
 J_1 &:=&   \int_{D_t \cap B(y(t),\lambda)} |x-y(t)|^2 |p_\tau(x,y(t))|^2 \rmd x ,\\
 J_2 &:=&    \int_{D_t \sauf B(y(t),\lambda)} |x-y(t)|^2 |p_\tau(x,y(t))|^2 \rmd x,\\
 \tilde J &:=& \int_{D_t} |p_\tau(x,y(t))|^2 \rmd x.
\end{eqnarray*}
 We then have
\begin{eqnarray}
\label{est.J1}
 J_1 &\le &  \lambda^2  \tilde J = 4\kappa^{-4} \dtyt^2 \tilde J.
\end{eqnarray}
 On the one hand thanks to\refq{est.ptau} we have 
\begin{eqnarray*}
 J_2 &\le & |D_t| 4\pi \kappa^{-4} \exp(-2\kappa\lambda \tau) ,
\end{eqnarray*}
 and, on the other hand,
\begin{eqnarray*}
 J_2 &\le & 4\pi \kappa^{-4} \int_{|x-y(t)|>\lambda} \!\!\! e^{-2\kappa \tau |x-y(t)|} \rmd x
 \le 16\pi^2 \kappa^{-4} \int_{r>\lambda}  e^{-2\kappa \tau r} r^2 \rmd r \\
 &\le& 20\pi^2 \lambda^2 \kappa^{-5} \tau^{-1}  \exp(-2\kappa \tau \lambda)
 \le 20\pi^2 \lambda^3 \kappa^{-4}  \exp(-2\kappa \tau \lambda)\\
 &\le& 20\pi^2 \kappa^{-10} (2\dtyt)^3  \exp(-2\kappa \tau \lambda).
\end{eqnarray*}
 Here we used $\lambda\ge \kappa^{-1}\tau^{-1}$.
 Hence
\begin{eqnarray}
\label{est.J2}
 J_2 &\le & 4\pi \kappa^{-4}\min(40\pi\kappa^{-6}  \dtyt^3,  |D_t|)\exp(-2\kappa \tau \lambda).
\end{eqnarray}
 Let us fix $x_0\in \partial D_t$ such that $\dtyt=|x_0-y(t)|$. Then $B(y(t),2\dtyt)\supset B(x_0,\dtyt)$
 so, thanks to\refq{est.ptau} and to (H3b), we have
\begin{eqnarray*}
 \tilde J &\ge &  \int_{D_t\cap B(y(t),2\dtyt)} |p_\tau(x,y(t))|^2 \rmd x \\
 &\ge & |D_t\cap B(y(t),2\dtyt)|  4\pi \kappa^4 (2\dtyt)^{-2} \exp(-4\kappa^{-1}\tau \dtyt)\\
 &\ge & \pi \kappa^4 |D_t\cap B(x_0,\dtyt)| \dtyt^{-2} \exp(-2\kappa\tau \lambda)\\
 &\ge & \pi \kappa^4  L_D \min(|D_t|,|B(x_0,\dtyt)|) \dtyt^{-2}  \exp(-2\kappa\tau \lambda)\\
 &\ge & \pi \kappa^4  L_D\min(|D_t|,\frac43\pi \dtyt^3)\dtyt^{-2}  \exp(-2\kappa\tau \lambda).
\end{eqnarray*}
 Then
\begin{equation}
\label{est.J2/tildeJ}
 \frac{J_2}{\dtyt^2\tilde J} \le C L_D^{-1} \kappa^{-14} ,
\end{equation}
 for some numerical parameter $C>0$.
 From\refq{est.J2/tildeJ} and\refq{est.J1} we obtain
$$ J\le C'L_D^{-1} \dtyt^2 \tilde J, $$
 which is the estimate to prove. 
 \qed

\medskip
 \ssl {Proof of Lemma \ref{lem2.maj-dqtau}}.
 It is the direct consequence of Lemma \ref{lem.maj-dqtau} and Lemma \ref{lem.(x-y)ptau}.

\medskip
\ssl {Proof of Lemma \ref{lem.est1-lnIinfty}}\\
 1) Thanks to\refq{est.ptau} and to Lemma \ref{lem.R9}  we have, for all $\tau>\tau_0$,
\begin{eqnarray*}
 I_\infty(\tau) &\le& c \tau^{-4}   \int_0^T \int_{D_t} 16 \kappa^{-4} \epsl^{-2}_\Sigma e^{-2\kappa\tau\epsl_\Sigma} \rmd x \rmd t + C_1 e^{-4\kappa^{-1} d_\Om\tau} \\
 &\le& C_2 \epsl_\Sigma^{-2}\tau^{-4} e^{-2\kappa\epsl_\Sigma\tau}.
\end{eqnarray*}
 We then obtain\refq{est1.lnIinfty}.
 
2)  Let us  fix $T_1 \in (0, \min(\alpha|D_\teta|^{1/3},\frac12\teta,\frac12(T-\teta),\alpha^2,
 \frac{\alpha^3}{\sqrt{8C_M}}))$ sufficiently small such that we have, thanks to (H2), 
\begin{equation}
 \label{min.|Dt|}
 |D_t|\ge \frac12 |D_\teta| >0 \quad  {\rm for} \quad |t-\teta|\le T_1.
\end{equation}
 Thanks to\refq{minIinfty01} in Lemma \ref{lem.R9} we have, for $\tau>\tau_0$,
$$
 I_\infty(\tau) \ge I_0  - R_0
$$
 with
\begin{eqnarray*}
 R_0 &:=& \frac1{c} \tau^{-4} \int_{|t-\teta|\ge T_1} C_M\dtyt^2 e^{-2\tau\mu|t-\teta|}  \int_{D_t}
  p^2_\tau(x,y(t)) \rmd x \rmd t \\
  && + C_1 e^{-\tau \mu\min(\teta,T-\teta)} ,\\ 
 I_0 &:=& \frac1{c} \tau^{-4} \int_{|t-\teta|\le T_1} (1-C_M\dtyt^2)e^{-2\tau\mu|t-\teta|}
 \int_{D_t}  p^2_\tau(x,y(t)) \rmd x \rmd t.
\end{eqnarray*}
 Thanks to\refq{est.ptau} we have
\begin{eqnarray*}
 R_0 &\le& C_3(\Sigma,v_0) \tau^{-4}  e^{-2\tau\mu T_1}.
\end{eqnarray*}
 By observing that, thanks to\refq{cond1.Sigma}, $|t-\teta|\le T_1$ implies 
 $\dtyt\le \frac2{\alpha^3}T_1 \le \frac1{\sqrt{2C_M}}$.
 Thus, putting
\begin{eqnarray*}
 B(t) := \int_{D_t\cap B(y(t),2\dtyt)}  e^{-4\tau\kappa^{-1} \dtyt} \rmd x,\quad \epsl\le |t-\teta|\le 2\epsl,
\end{eqnarray*}
 and restricting $\epsl$ to the interval $(0,T_1/2)$, we obtain
%
%
%
\begin{eqnarray*}
  I_0  &\ge& \frac1{2c'} \sup_{0\le r\le T}\{\rmd(y(r),D_r)\}^{-2}  \tau^{-4} \int_{|t-\teta|\le T_1}
   e^{-2\tau\mu|t-\teta|} B(t) \rmd t\\
 &\ge& c(M,\alpha,\Om,T,\kappa)  \tau^{-4} \int_{\epsl\le |t-\teta|\le 2\epsl} e^{-2\tau\mu|t-\teta|} 
 B(t) \rmd t.
\end{eqnarray*}
 Let us give a lower bound for $B(t)$, $\epsl\le |t-\teta|\le 2\epsl$.
 We have, thanks to\refq{cond1.Sigma},
$$
 D_t\cap B(y(t),2\dtyt) \supset D_t\cap B(x(t),\dtyt)  \supset D_t\cap B(x(t),\frac{\epsl}{2\alpha}),
$$
 for some $x(t)\in \partial D_t$. 
 Thanks to (H3b) we have
$$ |D_t \cap B(x(t),\dtyt)| \ge  L_D \min(|D_t|,\frac{\pi\epsl^3}{6\alpha^3}) .$$
 Thus, thanks to\refq{min.|Dt|} and since $\epsl\le \frac12 T_1 \le \frac12 \alpha |D_\teta|^{1/3}$,
 we have
$$
 |D_t| \ge \frac12 |D_\teta| \ge  \frac{\pi\epsl^3}{6\alpha^3} \quad {\rm for} \quad  \epsl\le |t-\teta|\le 2\epsl,
$$
 and so
\begin{equation}
\label{maj1.DtcapB}
  |D_t\cap B(y(t),2\dtyt)| \ge L_D \frac{\pi\epsl^3}{6\alpha^3} , \quad  \epsl\le |t-\teta|\le 2\epsl.
\end{equation}
 Then, thanks to\refq{cond1.Sigma},
$$
 B(t) \ge L_D \frac{\pi\epsl^3}{6\alpha^3} e^{-4\tau\kappa^{-1} \dtyt} \ge
 L_D \frac{\pi\epsl^3}{6\alpha^3} e^{-8\tau\kappa^{-1} \alpha^{-3} \epsl} ,
  \quad  \epsl\le |t-\teta|\le 2\epsl.
$$
 Finally we obtain
\begin{equation}
\label{maj.I0}
  I_0  \ge C(\Sigma_\epsl)  \tau^{-4} e^{-(8\kappa^{-1}\alpha^{-3}+4\mu)\epsl \tau},
\end{equation}
 with $C(\Sigma_\epsl) >0$.
 Let us put $\epsl_1= \min(\frac12T_1,\frac{\mu T_1}{\kappa^{-1}\alpha^{-3}+4\mu})$.
 For $\epsl\in (0,\epsl_1]$ we then have
\begin{eqnarray*}
  R_0/I_0  &\le & C'(\Sigma) e^{-\tau\mu T_1},
\end{eqnarray*}
 and so $R_0/I_0\le \frac12$ for $\tau\ge \tau_0$ (eventually modified).
 Thus, for $\epsl\in (0,\epsl_1]$,
$$
 I_\infty(\tau) \ge  \frac12 C(\Sigma_\epsl)  \tau^{-4} e^{-(8\kappa^{-1}\alpha^{-3}+4\mu)\epsl \tau}
$$
 which implies\refq{est2.lnIinfty}.

\qed


\begin{thebibliography}{99}

\bibitem{ARO.BOU}, D. G. Aronson, \textit{Bounds for the fundamental solution of a parabolic equation},
 Bull. Amer. Math. Soc. \textbf{73}, Number 6 (1967), 890-896. 

\bibitem{CRI.STA}
 M. Di Cristo and S. Vessella, \textit{Stable determination of the discontinuous conductivity
 coefficient of a parabolic equation}, arXiv:0904 (2009).

\bibitem{DAI.PRO}
 Y. Daido, H. Kang and G. Nakamura, \textit{A probe method for the inverse boundary value
 problem of non-stationary heat equations}, Inverse Problems \textbf{23} (2007), 1787-1800.

\bibitem{ELA.UNI}
 A. Elayyan and V. Isakov, \textit{On uniqueness of the recovery of the discontinuous
  conductivity coefficient of a parabolic equation},
 SIAM. J. Math. Anal. \textbf{28} (1997), 49-59.

%

\bibitem{FAB.NEW}
 E.B. Fabes and D.W. Stroock, \textit{A new proof of Moser's parabolic Harnack inequality
  via the old ideas of Nash}, Arch. Rat. Mech. Anal. ,  \textbf{96} (1986) 327-338.

\bibitem{FRI.UNI} A. Friedman and V. Isakov, \textit{On the uniqueness in the inverse
 conductivity problem with one measurement}, Indiana Univ. Math. J. {\bf 38} (1989), 563-579.
 
\bibitem{GAI.INV1}
 P. Gaitan, H. Isozaki, O. Poisson, S. Siltanen and J. Tamminen, \textit{Inverse problems
 for time-dependent singular heat conductivities - One dimensional case}, SIAM Journal of
 Mathematical Analysis 45(3), pp. 1675-1690 (2013).

\bibitem{GAI.INV2}
 P. Gaitan, H. Isozaki, O. Poisson, S. Siltanen, J. Tamminen, \textit{Inverse problems
 for time-dependent singular heat conductivities. Multi dimensional case}, 
 Communications in Partial Differential Equations 40(5), pp. 837-877 (2014).

\bibitem{GAI.PRO}
 P. Gaitan, H. Isozaki, O. Poisson, S. Siltanen and J. Tamminen, 
\textit{Probing for inclusions in heat conductive bodies},
 Inverse Problems and Imaging \textbf{6} (2012), pp. 423--446.

\bibitem{IKE.ENC}
 M. Ikehata and M. Kawashita, \textit{The enclosure method for the heat equation},
 Inverse Problems \textbf{25} (2009), 075005.1 .

\bibitem{IKE.PRO}
 M. Ikehata, \textit{Probe method and a Carleman function}, Inverse Problems \textbf{23}
 (2007), pp.1871--1894, doi:10.1088/0266-5611/23/5/006 .

\bibitem{IKE.REC}
 M. Ikehata and M. Kawashita, \textit{On the reconstruction of inclusions in a heat
 conductive body from dynamical boundary data over a finite time interval},
 Inverse Problems \textbf{26} (2010), 095004.
 
\bibitem{IKE.SIZ}
 M. Ikehata, \textit{Size estimation of inclusion},
  J. Inv. Ill-Posed Problems, \textbf{6}, N$^o$ 2 (1998), 127-140.
 
\bibitem{ISA.INV} V. Isakov, \textit{Inverse Problems for partial differential equations},
 Appl.Math.Sci., {\bf 127}, Berlin Springer (1988).

\bibitem{ISA.REC} V. Isakov, K. Kim, G. Nakamura, \textit{Reconstruction of an unknown
 inclusion by thermography}, Ann. Scuola Norm. Sup. Pisa Cl. Sci.,
 \textbf{Vol. IX}, issue 4  (2010), p. 725-758.

\bibitem{KAW.UNI}
 H. Kawakami and M. Tsuchiya, \textit{Uniqueness in shape identification of
 a time-varying domain and related parabolic equations on non-cylindrical domains}
 Inverse Problems \textbf{26} (2010), 125007 (34pp),
 doi:10.1088/0266-5611/26/12/125007

\bibitem{LIO.NON} J.L. Lions and E. Magenes \textit{Non-homogeneous boundary value problems and applications II}  
 Berlin: Springer (1972).

\bibitem{NAS.CON} J. Nash, \textit{Continuity of solutions of a parabolic and elliptic equations},
 Amer. J. Math., {\bf 80} (1958), p. 931-954.
 
\bibitem{POI.REC}
 O. Poisson,  \textit{Recovering time-dependent inclusion in heat conductive bodies using
 a dynamical probe method}, Journal of Mathematical Analysis and Applications, 441 (2):  
 862-844 (2016).

\bibitem{VES.QUA} S. Vessella, \textit{Quantitative estimates of unique continuation for parabolic equations,
 determination of unknown time-varying boundaries and optimal stability estimates},
 Inverse Problems, {\bf 24/2} (2008), 023001.
 
\bibitem{WLO.PAR} J. Wloka, \textit{Partial differential equations}, London: Cambridge University Press
(1987).
\end{thebibliography}
\end{document}